\documentclass[11pt]{article}

\usepackage[margin=1in]{geometry}
\usepackage{amsmath,amssymb,amsthm,mathtools}
\usepackage{bm}
\usepackage{booktabs, tabularx}
\usepackage{graphicx}
\usepackage{placeins}
\usepackage{enumitem}
\usepackage{microtype}
\usepackage{setspace}
\usepackage{xcolor}
\usepackage{natbib}
\usepackage{hyperref}
\usepackage[nameinlink,capitalise]{cleveref}

\hypersetup{colorlinks=true,linkcolor=black,citecolor=black,urlcolor=black}
\newtheorem{theorem}{Theorem}
\newtheorem{proposition}{Proposition}
\newtheorem{lemma}{Lemma}
\newtheorem{corollary}{Corollary}
\newtheorem{assumption}{Assumption}
\newtheorem{definition}{Definition}
\newtheorem{remark}{Remark}

\crefname{theorem}{Theorem}{Theorems}
\Crefname{theorem}{Theorem}{Theorems}
\crefname{proposition}{Proposition}{Propositions}
\Crefname{proposition}{Proposition}{Propositions}
\crefname{lemma}{Lemma}{Lemmas}
\Crefname{lemma}{Lemma}{Lemmas}
\crefname{corollary}{Corollary}{Corollaries}
\Crefname{corollary}{Corollary}{Corollaries}
\crefname{assumption}{Assumption}{Assumptions}
\Crefname{assumption}{Assumption}{Assumptions}
\crefname{definition}{Definition}{Definitions}
\Crefname{definition}{Definition}{Definitions}
\crefname{remark}{Remark}{Remarks}
\Crefname{remark}{Remark}{Remarks}
\crefname{example}{Example}{Examples}
\Crefname{example}{Example}{Examples}

\newcommand{\cF}{\mathcal F}

\newcommand{\cC}{\mathcal C}

\newcommand{\cT}{\mathcal T}

\newcommand{\CAAR}{\mathrm{CAAR}}
\newcommand{\E}{\mathbb E}
\newcommand{\Pp}{\mathbb P}

\title{Congestion-Aware Adaptive Routing for LLM Classification}
\author{%
Yuan Guo\\
{\small School of Business, The George Washington University}\\[1em]
Stefanus Jasin\\
{\small Stephen M. Ross Business School, University of Michigan}\\[1em]
Chen-An Lin\\
{\small Mitchell E. Daniels, Jr. School of Business, Purdue University}
}
\date{}

\begin{document}
\maketitle

\begin{abstract}
Organizations use large language models (LLMs) to classify jobs that differ in uncertainty and the cost of errors. Querying
additional models may improve accuracy, but it also incurs fees and increases congestion at shared servers. We study congestion-aware adaptive routing (CAAR) policies that adapt model selection and stopping decisions to observed outputs while accounting for the steady-state congestion they generate. The objective is
to minimize the long-run average cost of misclassification penalties, query fees, and holding costs. %
We can represent CAAR policies as a mixture of decision trees, and derive exact holding costs for first-come-first-served servers with exponential service times and processor-sharing
servers with general service times. We formulate a convex optimization problem with unique optimal model loads. Congestion prices capture the marginal system-wide holding costs of additional queries and guide routing and stopping decisions through a Bayesian dynamic program. We characterize optimal stopping decisions and identify models whose information does not justify their cost. A column-generation algorithm computes optimal CAAR policies and provides an explicit optimality-gap bound. We also bound the additional cost of limiting cascade depth relative to full-depth CAAR. Finally, under a common query-depth limit, cascades with a predetermined model order can be arbitrarily more costly than optimal CAAR, even when the order is randomized and stopping decisions adapt to observed outputs. Numerical experiments distinguish the
benefits of adaptive model selection and congestion-aware design.
\end{abstract}


\section{Introduction}
\label{sec:introduction}
Organizations deploying large language models (LLMs) increasingly face a choice among models that differ in classification accuracy, usage fees, and processing capacity \citep{chen2023frugalgpt,ong2024routellm,huang2025thriftllm}.
Routing determines which model to query and whether
to seek additional evidence before making a decision.
A model cascade addresses this problem by querying
models sequentially, using earlier outputs to
determine whether further queries are needed.
When many jobs share the same inference servers,
these choices affect not only the cost and accuracy
of serving individual jobs but also the delays
experienced throughout the system. Routing more jobs
to a highly accurate model may improve the 
quality, but it can also increase query expenditures
and congestion at that model. Similarly, consulting
multiple models for a single job can improve the
final classification, but each additional query
incurs a fee and adds workload to the system.

We focus on binary classification tasks, in which each
job must be assigned one of two labels. For example,
an LLM may screen a user prompt for safety risks or
determine whether a generated response violates a
content policy \citep{inan2023llama}. In legal context, an LLM can assess whether a sentence overrules an earlier decision, a task used to evaluate LLM selection methods \citep{huang2025thriftllm}. In both applications, classification errors and delays in reaching a decision can be costly. 
This focus allows us to obtain explicit stopping rules and information-value bounds for congestion-aware adaptive routing.
Section~\ref{sec:conclusion} discusses extensions to more general inference tasks and loss functions.

Existing LLM studies examine both output-adaptive model selection
and routing with workload-dependent latency. Table~\ref{tab:llm-positioning-matrix} distinguishes these features across representative formulations. 

\begin{table}[ht]
\centering
\caption{Comparison of Selected LLM Routing Studies.}
\label{tab:llm-positioning-matrix}
\footnotesize
\setlength{\tabcolsep}{5pt}
\renewcommand{\arraystretch}{1.2}
\begin{tabularx}{\linewidth}{
@{}>{\raggedright\arraybackslash}p{0.23\linewidth}
*{3}{>{\raggedright\arraybackslash}X}@{}}
\toprule
\textbf{Response time in policy design}
& \textbf{Single-model assignment}
& \textbf{Prescribed order with adaptive continuation}
& \textbf{Output-adaptive model selection and stopping} \\
\midrule

\textbf{No explicit response-time component}
&
RouteLLM \citep{ong2024routellm}
&
FrugalGPT \citep{chen2023frugalgpt}\newline 
ThriftLLM \citep{huang2025thriftllm}\newline
\citet{belloni2026pandora}\newline
&
\citet{khalaf2026circuits}
\\
\midrule

\textbf{Exogenous response times or latency costs}
&
--
&
TREACLE \citep{zhang2024treacle}\newline
&
\citet{dekoninck2025unified}\newline
\citet{li2026sequential}
\\
\midrule

\textbf{Endogenous congestion}\newline
&
RouterWise \citep{kasnavieh2026routerwise}\newline
&
Cascadia \citep{jiang2026cascadia}\newline
&
\textbf{This paper}\newline
\\
\bottomrule
\end{tabularx}
\end{table}

RouteLLM selects a single model for each request, while FrugalGPT uses a prescribed cascade \citep{ong2024routellm,chen2023frugalgpt}. A prescribed order may depend on the initial request context, and decisions to continue can still depend on observed outputs. Other studies adapt both model selection and stopping to earlier outputs, with query costs or response-time distributions specified externally \citep{dekoninck2025unified,khalaf2026circuits,li2026sequential}. RouterWise and Cascadia account for routing-dependent workloads. The former uses single-model assignment, while the latter optimizes deployment and continuation along a prescribed cascade \citep{kasnavieh2026routerwise,jiang2026cascadia}.
We jointly optimize output-dependent model selection and stopping,
accounting for the steady-state congestion these decisions
generate across jobs.

This setting is particularly relevant for organizations that operate inference servers in-house or through dedicated hosting arrangements. A McKinsey survey of more than 700 technology leaders and senior developers found that 63\% reported use of open or partially open AI models, while 40\% of leaders preferred hosting models on their own infrastructure \citep{McKinseyOpenSource2025}. Similarly, BCG identifies on-premises and single-tenant cloud deployments as alternatives to shared cloud services \citep{BCGFactory2023}. These deployment settings motivate our focus on organizations that route inference requests across resources with specified service capacities, where routing decisions directly affect server workloads and congestion.

More specifically, we ask \emph{how an organization should design adaptive cascade routing, accounting for both the information gained from each query and the congestion it creates}. To study this question, we consider a system %
with multiple job types, each characterized by its own prior probabilities and misclassification penalty. Jobs arrive according to independent Poisson processes and are processed by shared multiple LLM servers. Our main model assumes first-come-first-served (FCFS) service for each model with exponential service times, but the results also extend to processor-sharing servers with general service times. Each queried model receives the same original job input and returns a noisy binary signal. Earlier outputs guide routing and stopping decisions but are not included in subsequent model inputs. This aligns with sequential classification workflows
such as assessing whether a response violates a safety policy or
whether a legal sentence overrules an earlier decision
\citep{khalaf2026circuits,huang2025thriftllm}.
A job may be classified immediately or query up to \(L\) distinct models before making a final classification.

We study congestion-aware adaptive routing (CAAR) for sequential classification. We optimize over policies whose routing and stopping decisions may be randomized and depend on a job's type, previously queried models, and observed outputs, but not on current queue states.
Our objective is to minimize the long-run average cost of misclassification penalties, query fees, and holding costs. Congestion enters through the steady-state model loads generated by routing decisions across all job types. Our results characterize and compute an optimal policy within this admissible class.

\paragraph{Contributions.}
Our first contribution is a system-level framework for designing
congestion-aware adaptive routing (CAAR). We jointly optimize routing and
stopping rules across job types to balance misclassification penalties,
query fees, and steady-state holding costs. Within the admissible policy class,
we derive exact holding costs under our queueing specifications. Although
a job's route depends on earlier outputs, these costs depend only on
the aggregate visit rate at each model. By representing policies as
mixtures of decision trees, model loads are linear in the mixture
probabilities. This yields a convex optimization problem with a unique
optimal load vector. We further characterize optimal policies through
congestion prices that capture the marginal system-wide holding cost
of additional queries. These prices separate policy selection across
job types. System optimality requires that the selected policies
generate loads consistent with the prices used to select them.

Our second contribution is to characterize and compute optimal
routing policies and quantify the cost of limiting their depth.
At given congestion prices, a Bayesian dynamic program identifies
a least-cost routing tree for each job type. We establish a
two-sided threshold structure for stopping. Under positive signal
orientation, information-value bounds rule out model queries whose
potential benefits cannot justify their fees and congestion costs. The same dynamic program supplies improving trees for a column-generation algorithm that jointly solves policy mixtures and congestion prices. This approach avoids enumerating all trees in advance. An explicit bound on the remaining cost improvement allows computation to terminate with a certified optimality gap within the policy class. We then bound the additional system cost of limiting query depth based on system primitives. This bound decreases with the permitted depth and identifies a sufficient depth for a prescribed cost tolerance
relative to full-depth CAAR. Together, these results allow shorter
cascades to be computed with explicit bounds on both the cost of
restricting depth and the remaining optimization error.

Our third contribution is to quantify the value of flexibility and congestion awareness in routing design. Comparisons with fixed and randomized model rankings examine the value of choosing subsequent models from observed outputs. We establish a worst-case separation: committing to a model order in advance can be arbitrarily more costly than CAAR, even when stopping adapts to observed outputs. Separately, comparisons with policies that ignore load-dependent waiting costs examine the value of congestion-aware design. In our two-type specialist experiments, CAAR reduces total cost by 5.80\% on average and up to 24.42\% relative to the
service-time-only benchmark. We also identify settings in which
simpler policies perform nearly as well, clarifying when additional
routing flexibility provides meaningful benefits.

Finally, Appendix~\ref{apdx:repeated-queries} extends the framework
to repeated queries of the same model. Each call rejoins the model's
queue as a separate request. The queueing, convex optimization, and
congestion-price results continue to hold. %
The Bayesian formulation also extends when repeated outputs are conditionally independent given the true label and have the same likelihoods across calls.

\paragraph{Organization of the paper.}
The remainder of the paper is organized as follows.
Section~\ref{sec:literature} reviews the related
literature, and Section~\ref{sec:policy-hierarchy}
introduces the model and CAAR policies.
Section~\ref{sec:system-foundation} establishes the
exact queueing costs, convex optimization formulation,
and congestion-price characterization.
Section~\ref{sec:policy-structure} studies optimal
routing and stopping decisions.
Section~\ref{sec:column-generation} present the column generation procedure, and Section~\ref{sec:effective-depth} bounds the cost
of limiting query depth.
Section~\ref{sec:numerical} compares CAAR with
fixed-ranking cascade and congestion-blind policies through analytical and
numerical results.
Section~\ref{sec:conclusion} concludes the paper.
All proofs are provided in the appendix. %

\section{Literature Review}
\label{sec:literature}

Our work draws on three streams of literature:
LLM routing and cascading, sequential information
acquisition, and queueing networks. We review
each in turn, highlighting how our work brings
these streams together by studying routing and
information acquisition decisions when the cost
of obtaining additional information depends on
congestion at shared servers.

\subsection{LLM routing and adaptive inference}

Beyond the routing studies summarized in Table~\ref{tab:llm-positioning-matrix}, related work examines how additional model outputs are combined or revised, when to stop sampling, and how much deeper inference can improve performance. LLM-Blender combines responses from multiple models
\citep{jiang2023llm}.
Self-consistency aggregates sampled answers, whereas self-refinement
revises an answer through feedback
\citep{wang2022selfconsistency,madaan2023selfrefine}.
\citet{huang2026stopping} use Bayesian stopping to identify the most
consistent answer, while BEACON balances expected reward improvements
against sampling costs \citep{wan2025beacon}.
We instead acquire classification evidence about an unchanged job
input, whose true label need not be the most frequent model answer.
\citet{bouchard2026depth} relates limited empirical gains from deeper
cascades to the information in confidence scores. Our effective-depth
bound concerns the total system cost of truncating a CAAR policy,
including its effect on congestion.

\subsection{Sequential information acquisition
and gatekeeping}

A central question in sequential information
acquisition is whether the benefit of obtaining
another signal justifies its cost. This trade-off
has been extensively studied in sequential testing
and Bayesian decision making
\citep{Wald1947,DeGroot1970,Chernoff1959}.
Cost-sensitive classification and budgeted feature
acquisition apply this trade-off to information
sources with fixed acquisition costs
\citep{trapeznikov2013supervised,nan2017adaptive}.
The literature on finite-horizon partially
observed control also provides the mathematical
foundation for the concavity and piecewise-affine
value functions that arise in our analysis
\citep{Astrom1965,SmallwoodSondik1973}.

Related questions have also been studied in
service operations, particularly in diagnostic
systems and gatekeeper networks. 
\citet{alizamir2013diagnostic} examine the
trade-off between diagnostic accuracy and
congestion, while gatekeeper models study referrals
from a primary resource to a specialist
\citep{shumsky2003gatekeepers,hasija2005staffing}.
Empirical work documents workload effects on testing,
task ordering, and referrals
\citep{kc2009impact,batt2017early,
ibanez2018discretionary,freeman2021gatekeeping}.
We study a network of heterogeneous information
sources in which a job can collect several signals
before classification. Congestion enters through
the policy's steady-state traffic rather than
observed queue lengths at each decision.

\subsection{Queueing networks, routing,
and congestion prices}

Classical routing policies such as
join-the-shortest-queue use real-time queue lengths
to direct jobs toward less congested servers
\citep{winston1977optimality,ephremides1980routing}.
Recent LLM research also connects routing and learning
with queueing.
\citet{bae2026queueing} study LLM routing and scheduling
with user retrials, learning preferences from
feedback while maintaining queue stability. 
\citet{ao2026service} study congestion caused by
KV-cache growth in memory-constrained LLM serving.
They examine admission control and request mixing
to mitigate eviction-induced throughput losses.
Our focus is on model outputs as evidence about a
job's true label. Models are not interchangeable
servers, and a job may consult several of them.
Routing therefore determines both the information
collected for classification and the resulting
service workload.

Our queueing analysis builds on the theory of
multiclass product-form networks
\citep{BCMP1975,Kelly1979}. 
We also use marginal holding costs as congestion
prices, following network-routing and flow-control
principles \citep{Kleinrock1964,BertsekasGallager1992}.
Our contribution brings these service-system
foundations together with evidence-dependent querying
and stopping. %

\section{Model, CAAR Policies, and Objective}
\label{sec:policy-hierarchy}
We consider an organization that processes binary classification jobs using multiple LLMs on servers with specified service capacities. These servers are operated in-house or through dedicated hosting arrangements. The organization's routing decisions determine their
workloads and resulting congestion. Our framework does not directly
cover shared external APIs whose delays also depend on substantial
traffic from other users.

For each job, the router decides whether to classify it immediately or query one or more models.
After each query, it can either stop and make a final
classification or consult another model. These
decisions involve a trade-off between classification
accuracy, query fees, and congestion-induced delays.
We first describe the job characteristics, model
outputs, and service system, and then formally define
the admissible routing policies and the objective
function.

\subsection{Jobs, service, and Bayesian signals}

There are $J\ge1$ job types, indexed by
$j\in[J]:=\{1,\ldots,J\}$, and $M\ge1$ LLMs,
indexed by $m\in[M]:=\{1,\ldots,M\}$.
Jobs of type $j$ arrive according to a Poisson
process with rate $\lambda_j>0$, independently
across job types. The total arrival rate is
$\lambda:=\sum_{j=1}^J\lambda_j$.
Each job's type is observed upon arrival, but
its true label is unknown. For a generic type-$j$
job, let $Y_j\in\{A,B\}$ denote its true label,
with prior probability
$\Pp(Y_j=A)=\rho_j\in(0,1)$.
For example, the two labels may indicate whether
a prompt is safe or unsafe. %
An incorrect
classification of a type-$j$ job incurs a penalty
$b_j>0$. Thus, job types may differ in both
their prior probabilities and the costs of
misclassification. 

Each LLM operates as a work-conserving,
first-come-first-served (FCFS) single-server queue. 
Service times at model $m$ are independent and
exponentially distributed with rate $\mu_m>0$,
regardless of the job type. (We use FCFS as our
baseline, but the results
also extend to processor-sharing servers with general
service times, as discussed in Remark \ref{remark1} in Section \ref{sec:system-foundation}.) Each query to model $m$
incurs a fee $c_m>0$, and the model's output becomes
available upon service completion. \textcolor{black}{For an in-house deployment, this fee is interpreted as a marginal operating cost rather than a payment to an external provider.} In addition,
all jobs incur a common holding cost at rate $h>0$
per unit of time spent in the system, including
both waiting and service. Models then differ in
their query fees, processing capacities, and
the quality of the information they provide described below.

For a generic type-$j$ job, let
$X_{jm}\in\{A,B\}$ denote the output that model
$m$ would produce when queried with the job’s original input. Outputs from earlier queries are available to the router but are not included in the input to model m. The router observes $X_{jm}$ upon completion of the query.  %
We characterize the model's accuracy separately
for the two possible true labels:
\begin{align*}
q_{jm}^A := \Pp(X_{jm}=A\mid Y_j=A), \quad
q_{jm}^B
= \Pp(X_{jm}=B\mid Y_j=B).
\end{align*}
Thus, $q_{jm}^A$ and $q_{jm}^B$ represent the
probabilities that model $m$ correctly classifies
a type-$j$ job when its true label is $A$ and
$B$, respectively. These accuracies may differ
across labels, models, and job types, allowing
the same model to be more useful for some types of 
jobs than others.

\begin{assumption}[Interior independent signals and independent service]
\label{ass:signals}
For every $j,m$, $q_{jm}^A,q_{jm}^B\in(0,1)$. Conditional on $Y_j$, the outputs $(X_{j1},\ldots,X_{jM})$ are mutually independent. Labels and potential outputs are independent across jobs. Arrival processes, service times, and policy randomization are mutually independent and independent of all labels and potential outputs. A model's output law does not depend on the time or position at which it is queried.
\end{assumption}

Conditional independence specifies how evidence from separate model queries can be combined. Given the true label, each unused
model retains its specified output distribution after earlier outputs are observed. The router can therefore combine evidence by multiplying likelihoods and use its current belief about the label to predict the output of an unused model. It still permits unconditional dependence among outputs through their shared dependence on the true label. Independence from service times further separates information quality from
processing duration. Service times contribute to delay costs but provide no additional evidence about the job's true label.

For each job type $j$ and model $m$, define
the conditional likelihood of output $x$ given
the true label $y$ as
$\ell_{jm}^y(x):=\Pp(X_{jm}=x\mid Y_j=y)$.
Using the conditional accuracies defined above,
these likelihoods are
$\ell_{jm}^A(A):=q_{jm}^A$, 
$\ell_{jm}^A(B):=1-q_{jm}^A$, 
$\ell_{jm}^B(A):=1-q_{jm}^B$, and
$\ell_{jm}^B(B):=q_{jm}^B$.
Suppose a type-$j$ job has already queried some
models and observed their outputs. Let
$p\in[0,1]$ denote the resulting posterior
probability that its true label is $A$.
If the job next queries an unused model $m$,
we denote by $s_{jm,x}(p)$ the probability
that the model produces output $x$, and by
$\tau_{jm,x}(p)$ the updated posterior
probability of label $A$ after observing
that output. By the law of total probability
and Bayes' rule, respectively,
\begin{align}
s_{jm,x}(p)
&:=p\ell_{jm}^A(x)+(1-p)\ell_{jm}^B(x),
\label{eq:predictive}
\\
\tau_{jm,x}(p)
&:=\frac{p\ell_{jm}^A(x)}{s_{jm,x}(p)}.
\label{eq:update}
\end{align}
Under Assumption~\ref{ass:signals}, the
denominator in \eqref{eq:update} is strictly
positive, so the posterior update is
well-defined for every $p\in[0,1]$.

\subsection{Congestion-aware adaptive routing (CAAR)}
A job queries models sequentially, and one at a time.
Upon arrival, the router may either classify the
job immediately or send it to a model for an
initial assessment. After each query is completed,
the router observes the model's output and decides
whether to stop and make a final classification
or query another model. We assume that routing
and classification decisions take no additional
time. We impose a maximum query depth
$L\in\{1,\ldots,M\}$, meaning that each job may
query at most $L$ distinct models, with each
model queried at most once. A job must stop
when this limit is reached, although it may
stop earlier. We also use $L=0$ to represent
the case in which every job is classified
immediately without querying any model.

An important feature of our framework is that
routing decisions can adapt to the information
collected from previous queries. For example,
two jobs of the same type may be sent to different
models after receiving different outputs from
their initial queries. Similarly, one job may
have enough evidence to stop, while another may
benefit from further queries.

We focus on congestion-aware adaptive routing (CAAR) policies,
which adapt routing and stopping decisions to a job's
observed history without using real-time queue states.
The routing rules are chosen to account for the
steady-state congestion they generate across jobs.
We formally define this policy class below.

\begin{definition}[{\color{black}Admissible policies for congestion-aware adaptive routing}]
\label{def:sa-class}
For a given depth limit $L$, the class $\Pi_L^{\CAAR}$ consists of routing, stopping, and terminal classification rules that depend only on a job's observed type and the history of models it has queried and their outputs. Decisions may be randomized, with randomization probabilities depending on the observed history. All jobs of the same type follow the same rule, using private randomization that is independent across jobs and independent of labels, model outputs, arrivals, and service times. In particular, decisions do not %
depend on current queue lengths, queue composition, FCFS manner, or elapsed waiting or service times.
\end{definition}
\textcolor{black}{The definition specifies the information a policy may use during execution. We select a policy from this class using the system's arrival rates, service capacities, query fees, and signal distributions.
Congestion awareness enters through the objective, which includes the steady-state holding costs generated by routing decisions across jobs.
Policies optimized under other objectives can also belong to this class, as illustrated by the benchmarks in Section~\ref{sec:numerical}. }
A fully queue-state-dependent policy could also respond to the evolving service system. Optimizing over that broader class would generally require tracking unfinished jobs' information states together with their positions in the FCFS queues. Our optimization instead focuses on the CAAR class. In Appendix~\ref{sec:repeated-queries}, we extend this framework to allow repeated queries to the same model, with each new call joining the model's queue as a separate request.

\subsection{The cost objective}
\label{subsec:CAAR-objective}

Our objective is to find a CAAR policy
$\pi\in\Pi_L^{\CAAR}$ that minimizes the
long-run average cost of classification errors,
query fees, and holding costs. For each job
type $j$, let $e_j(\pi)$ denote the probability
of a classification error under policy $\pi$,
and let $v_{jm}(\pi)$ denote the probability
that a type-$j$ job queries model $m$.
Since each model can be queried at most once,
$v_{jm}(\pi)$ also represents the expected
number of visits to model $m$ by a type-$j$
job. Aggregating across job types, the nominal
arrival rate of queries to model $m$ is
$\Lambda_m(\pi):=\sum_{j=1}^J\lambda_jv_{jm}(\pi)$, we also refer to $\Lambda_m$ as model $m$'s load.
Under stability, this quantity equals the
long-run average visit rate at model $m$.
The routing policy and signal distributions
jointly determine these visit rates, while
the congestion and delays at each model depend
on the total traffic it receives from all
job types.

For a policy $\pi$ whose induced queueing network
admits a stationary distribution, let
$\mathcal N_m$ denote the stationary number
of jobs at model $m$, including those waiting
in the queue and those receiving service.
Let $\mathbb E_\pi$ denote expectation under
this stationary distribution. The total
long-run average cost per unit time is
\begin{equation}
C_L(\pi)
:=\sum_{j=1}^J\lambda_j
\left[b_je_j(\pi)+\sum_{m=1}^Mc_mv_{jm}(\pi)\right]
+h\sum_{m=1}^M\mathbb E_\pi[\mathcal N_m].
\label{eq:policy-cost}
\end{equation}
The first term represents the expected
misclassification penalties and query fees
incurred per unit time across all job types.
The second term captures congestion costs %
as the holding cost rate $h$ times
the expected number of jobs in the system.
We define $C_L(\pi)=+\infty$ for any policy
that induces an unstable queueing network
or infinite expected stationary occupancy.
Our objective is to find a CAAR policy that
minimizes this total cost rate:
\begin{equation}
C_L^*
:=\inf_{\pi\in\Pi_L^{\CAAR}}C_L(\pi).
\label{eq:CAAR-optimal-value}
\end{equation}

Because immediate classification is allowed, a
feasible policy with finite cost always exists.
In particular, classifying each job using the
more likely label under its prior generates no
query traffic and incurs a total cost rate of
$\sum_{j=1}^J\lambda_jb_j\min\{\rho_j,1-\rho_j\}$.
Thus, we do not need to impose an aggregate
capacity condition such as
$\lambda<\sum_m\mu_m$ to ensure feasibility.
Instead, stability depends on the query traffic
generated by the chosen routing policy.
Furthermore, increasing the depth limit $L$
expands the set of feasible policies, since
jobs can always stop before reaching the limit.
Consequently, the optimal cost $C_L^*$ is
non-increasing in $L$.

To evaluate the objective in
\eqref{eq:policy-cost}, we need to determine
the expected number of jobs at each model.
This is not immediately straightforward because
the models visited by a job may depend on the
outputs of earlier queries. In the next section,
we show that the steady-state holding cost can
nevertheless be characterized and allows us to reformulate \eqref{eq:CAAR-optimal-value} as an explicit convex optimization problem. %

\section{Queueing Representation and Congestion-Aware Optimization}
\label{sec:system-foundation}

We now turn to evaluating and optimizing the total
cost of a CAAR policy in \eqref{eq:CAAR-optimal-value}.
The main challenge is that the sequence of models
queried by a job depends on previously observed
outputs, while the expected delay at each model depends
on the traffic generated by all jobs in the
system. To address this challenge, we represent
each routing rule as a decision tree that specifies
what action to take after every possible history
of model outputs. We then show that, despite
this output-dependent routing, the steady-state
holding cost depends only on the aggregate visit
rates at individual models. By allowing jobs to
be randomly assigned to different decision trees,
we obtain a convex optimization formulation.
We further characterize its optimality conditions
using congestion prices, which capture
the system-wide delay costs of additional queries
and allow routing decisions to be optimized
separately for each job type. 

\subsection{Representing CAAR policies
as tree mixtures}

A deterministic CAAR policy can be viewed as a complete adaptive routing plan as a decision tree that specifies what to do after every possible sequence of queried models and observed outputs. For example, a plan may first query model 1, stop with label $A$ if its output is $A$, and query model 2 if its output is $B$. The plan also specifies the final decision after model 2. Such a plan can be represented as a decision tree, where each node specifies an action and each observed output determines which branch the job follows. The outcome remains adaptive because different outputs can lead to different actions.

To formalize this representation, let
$\mathfrak h_t :=
((m_1,x_1),\ldots,(m_t,x_t))$
denote the history of a job after $t$ queries,
where $m_s$ is the model queried at step $s$
and $x_s\in\{A,B\}$ is its observed output.
Because each model can be queried at most once,
the models $m_1,\ldots,m_t$ are distinct.
We write $\mathfrak h_0:=\varnothing$ for
the initial history, before any query has
been made. We now formally define a
deterministic routing 
tree. We will use the terms decision tree and routing tree interchangeably. %

\begin{definition}[Tree representation]
\label{def:tree}
A complete deterministic depth-$L$ routing tree
assigns an action to every feasible history.
At a history of depth $t<L$, the action is
either to stop and assign a terminal label
in $\{A,B\}$ or to query a model outside
$\{m_1,\ldots,m_t\}$. At depth $L$, the job
must stop and be assigned a terminal label.
The root action is specified at
$\mathfrak h_0=\varnothing$. Let $\cT_L$
denote the set of all such trees, which
is finite because the number of models
and the maximum query depth are finite.
\end{definition}

A randomized CAAR policy instead uses independent private randomness to choose among feasible actions, with probabilities that may depend on the observed history. Randomization can occur as the job proceeds through the system. Alternatively, the job can randomly select a complete deterministic decision tree at arrival and follow it until departure. The selected tree fixes the action at every possible history, but the observed outputs determine which route the job would follow. The following lemma establishes the equivalence of these implementations. %

\begin{lemma}[Randomized admissible policies are exactly tree mixtures]
\label{lem:tree-mixture}
For every policy in $\Pi_L^{\CAAR}$ and every
job type, there is a probability distribution
over $\cT_L$ with the same joint law of queried
models, stopping histories, and terminal
decisions. Conversely, every distribution
over $\cT_L$ is implementable by a policy
in $\Pi_L^{\CAAR}$.
\end{lemma}

By this equivalence, we can represent any
randomized CAAR policy as a mixture of
deterministic routing trees. For each job
type $j$, let $z_{jT}$ denote the probability
of selecting tree $T\in\cT_L$, where
$z_{jT}\ge0$ and
$\sum_{T\in\cT_L}z_{jT}=1$.
Thus, designing a CAAR policy is equivalent
to choosing these probabilities for each
job type.

\begin{figure}[htbp]
\centering
\includegraphics[width=0.8\textwidth]{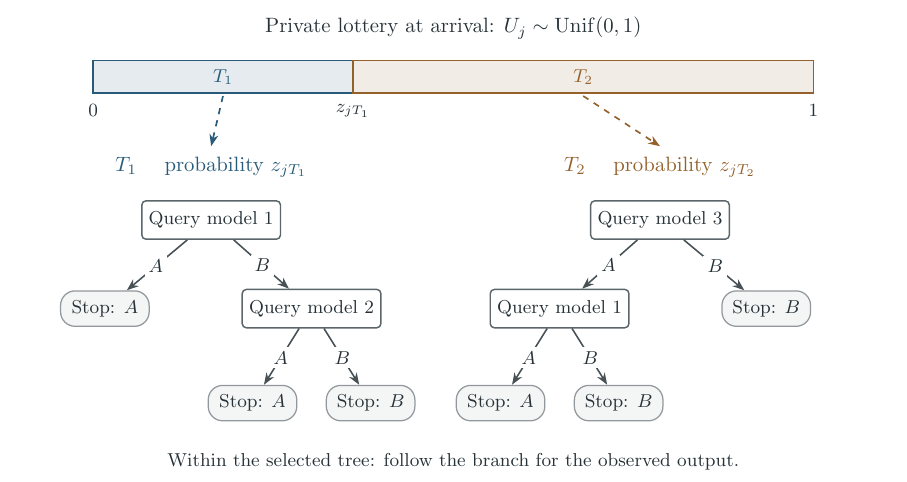}
\caption{An illustrative type-$j$ mixture with only two positive tree weights; all other trees have zero weight. A type-$j$ job draws $T_1$ with probability $z_{jT_1}$ and $T_2$ with probability $z_{jT_2}=1-z_{jT_1}$.}%
\label{fig:tree-mixture}
\end{figure}

To implement the mixture, partition $[0,1]$
into intervals of lengths $z_{jT}$, with
each interval corresponding to a tree.
Upon arrival, a type-$j$ job draws
$U_j$ uniformly from $[0,1]$ and follows
the tree associated with the interval
containing $U_j$. These random draws
are independent across jobs and independent
of labels, model outputs, arrivals, and
service times. Figure~\ref{fig:tree-mixture}
illustrates this construction for a mixture
of two trees.

Although different job types may use the same
routing tree, their error probabilities and model
visit probabilities can differ because of their
priors and signal distributions. For each job
type $j$ and tree $T$, define
\begin{align*}
e_{jT}&:=\mathbb P(\text{tree }T\text{ makes an error for type }j),
\\
v_{jTm}&:=\mathbb P(\text{tree }T\text{ queries model }m
\text{ for type }j),
\\
\omega_{jT}&:=b_je_{jT}+\sum_{m=1}^Mc_mv_{jTm}.
\end{align*}
Here, $\omega_{jT}$ is the expected cost of
following tree $T$, including the classification
error penalty and query fees, but excluding
holding costs. Since each model can be queried
at most once, $v_{jTm}\in[0,1]$ also represents
the expected number of visits to model $m$.
Moreover, $\sum_m v_{jTm}\le L$ because each
job makes at most $L$ queries. For any given
job type and tree, these quantities can be
calculated directly from the prior probabilities,
signal distributions, error penalties, and
query fees.

We can now express the traffic generated by a
tree mixture. Recall that $z_{jT}$ is the
probability that a type-$j$ job selects tree
$T$. The resulting aggregate load at model $m$
is
\begin{equation}
\Lambda_m(\bm z)
:=\sum_{j=1}^J\lambda_j\sum_{T\in\cT_L}z_{jT}v_{jTm}.
\label{eq:load-tree}
\end{equation}
To understand this expression, note that
type-$j$ jobs arrive at rate $\lambda_j$,
select tree $T$ with probability $z_{jT}$,
and visit model $m$ under that tree with
probability $v_{jTm}$. Summing the resulting
visit rates across all trees and job types
gives $\Lambda_m(\bm z)$. In particular,
the model loads are linear in the mixture
probabilities $\bm z$. Under stability,
$\Lambda_m(\bm z)$ equals the long-run
average visit rate at model $m$.

\subsection{Exact queueing costs under output-dependent routing}

We next evaluate the steady-state holding costs generated by the CAAR  policy. Although the model visit rates describe aggregate traffic, they do not by themselves determine expected queue occupancies. We therefore construct a queueing-network representation that accounts for each job's output-dependent routing and stopping decisions. Under our model assumptions, this representation yields exact holding costs as functions of the model visit rates.

Fix a tree mixture $\bm z$. For the purpose
of analysis, imagine that each arriving job
draws its routing tree, true label, and all
potential model outputs at the time of arrival.
Together, these draws determine the complete
route the job will follow, either departing
immediately without a query or visiting a
sequence $(m_1,\ldots,m_d)$ of distinct models,
where $1\le d\le L$. Importantly, this is
only an analytical construction. In the actual
system, the router still observes each model's
output only after its query is completed.
Because CAAR decisions do not depend
on queue states or service times, and the
draws are independent across jobs and of the
arrival and service processes, determining
the route in advance does not change the
system's stochastic behavior.

We refer to the information drawn for each
job as its \emph{arrival mark}. Since these
marks are independent across jobs and of
the arrival processes, the Poisson arrivals
of each job type can be split into independent
Poisson streams according to the routes
determined by their marks. Thus, although
routing decisions are adaptive in the actual
system, we can analyze the resulting traffic
as if each job were assigned a fixed route
at arrival.

To track jobs along these routes, consider
a type-$j$ job assigned to route
$(m_1,\ldots,m_d)$. At its $t$-th visit,
we assign it class
$(j,(m_1,\ldots,m_d),t)$.
The job receives service at model $m_t$,
then advances to the next stage or departs
after its final visit. These classes
distinguish jobs by their type, route,
and current stage, but do not change how
they are served. At each model $m$, all
classes share the same FCFS server with
exponential service rate $\mu_m$.
The resulting multiclass network therefore
satisfies the conditions of the open BCMP
theorem \citep{BCMP1975,Kelly1979}, which
allows us to characterize its stationary
queue lengths and holding costs.

\begin{proposition}[Exact open-BCMP representation]
\label{prop:bcmp-general}
Fix a tree mixture $\bm z$ and use the route-stage class expansion defined above. The resulting network is an open BCMP network whose physical model nodes are type-1 FCFS stations with class-independent exponential service rates. If $\Lambda_m(\bm z)<\mu_m$ for every $m\in[M]$, then the expanded network is positive recurrent and has the BCMP product-form stationary distribution. The stationary number of jobs waiting or in service at model $m$, denoted by $\mathcal N_m$, is geometric with mean
\begin{equation*}
\mathbb E[\mathcal N_m]
=\frac{\Lambda_m(\bm z)}{\mu_m-\Lambda_m(\bm z)}.
\end{equation*}
Define the corresponding holding-cost function by
\begin{equation*}
\Psi_m(\Lambda):=\frac{h\Lambda}{\mu_m-\Lambda},
\qquad 0\le\Lambda<\mu_m.
\end{equation*}
The total steady-state holding-cost rate is
$\sum_{m=1}^M\Psi_m(\Lambda_m(\bm z))$.
\end{proposition}

The expected number of jobs at each model therefore takes the familiar $M/M/1$ form. Importantly, Proposition~\ref{prop:bcmp-general} establishes this result through the product-form analysis of the entire multiclass network, without separately assuming that arrivals to each individual model form a Poisson process.

For a fixed CAAR policy, each job's route is determined by its selected tree and potential outputs. Under Assumption~\ref{ass:signals}, these can be sampled independently of service requirements and the queue state just before arrival. This property enables the preceding queueing representation while preserving output-dependent routing. Consequently, aggregate model loads suffice to evaluate steady-state holding costs under our assumptions. For the subsequent
optimization, we define $\Psi_m(\Lambda)=+\infty$
whenever $\Lambda\ge\mu_m$, consistent with
our convention that unstable policies incur
infinite cost. While FCFS provides a tractable baseline for modeling congestion, Remark 1 extends the queueing analysis to processor sharing, which may better capture LLM inference systems in which multiple requests concurrently share inference capacity.

\textcolor{black}{
\begin{remark} \label{remark1}
    While our main model centers on FCFS systems, our analysis and discussions can be extended to the processor-sharing system. The same holding-cost formula applies under equal processor sharing. In this variant, each request requires an independent amount of full-speed service to the model $m$ drawn from a general distribution with mean $1/\mu_m$, and simultaneous requests share the model’s capacity equally. The route-stage network then has BCMP processor-sharing stations \citep{BCMP1975}. Its stationary model occupancies yield the same function $\Psi_m$, so the optimization and subsequent congestion-price results can be extended. %
\end{remark}
}
\subsection{The convex system problem and unique optimal loads}

We can now formulate the CAAR problem
in terms of the tree-mixture probabilities.
Combining Lemma~\ref{lem:tree-mixture} and
Proposition~\ref{prop:bcmp-general}, we obtain
the following equivalent formulation of
\eqref{eq:CAAR-optimal-value}, which we call
the \emph{master problem}:
\begin{equation}
\label{prob:master}
\begin{aligned}
\min_{\bm z}\quad
C_L(\bm z) :=\;&
\sum_{j=1}^J\lambda_j
\sum_{T\in\cT_L}z_{jT}\omega_{jT}
+\sum_{m=1}^M\Psi_m(\Lambda_m(\bm z))
\\
\text{s.t.}\quad
&\sum_{T\in\cT_L}z_{jT}=1,
\qquad j\in[J],
\\
&z_{jT}\ge0,
\qquad j\in[J],\ T\in\cT_L,
\end{aligned}
\end{equation}
Here, $z_{jT}$ represents the probability that
a type-$j$ job follows tree $T$. The first term
of the objective captures the expected rate
of misclassification penalties and query fees,
while the second represents the steady-state
holding costs at the model loads defined in
\eqref{eq:load-tree}. The constraints ensure
that the mixture probabilities form a valid
distribution over trees for each job type.
Because jobs may be classified immediately
without querying any model, a finite-cost
feasible solution always exists. The convention that $\Psi_m(\Lambda)=+\infty$
for $\Lambda\ge\mu_m$ ensures that unstable
mixtures cannot be optimal.

An important property of the master problem
is its convexity. For each job type $j$ and
tree $T$, the coefficients $\omega_{jT}$ and
$v_{jTm}$ are fixed and do not depend on the
mixture probabilities. Consequently, both
the error-and-fee component of the objective
and the model loads are linear in $\bm z$.
On the other hand, the holding-cost functions $\Psi_m$ are strictly convex in their respective
model loads $\Lambda_m$ and thus convex in $\bm z$. These properties lead to the
following result.

\begin{theorem}[Convexity, existence, and unique optimal loads]
\label{thm:convexity}
The master problem \eqref{prob:master} is
convex in $\bm{z}$ %
and admits an optimal solution, with
optimal value $C_L^*$. The optimal tree
mixture $\bm z^*$ need not be unique, but the induced
optimal load vector
$\bm\Lambda^*=(\Lambda_1^*,\ldots,\Lambda_M^*)$
is unique.
\end{theorem}

The theorem establishes that although
different mixtures of routing trees may
achieve the same minimum cost, they must
generate exactly the same aggregate traffic
at each model. The optimal allocation
of model loads is then uniquely determined, even
when the implemented routing policies are not. Furthermore, every optimal model load is strictly below its service capacity.
This follows because holding costs diverge as model loads approach capacity, whereas classifying all jobs immediately provides a feasible policy with finite cost.

\subsection{Congestion prices and decomposition across job types}

We next characterize which routing trees are
used in an optimal solution to the master
problem \eqref{prob:master}. Notice that jobs of different types share the
same servers, so routing more traffic to one
model increases the delays experienced by
other jobs. To capture this interaction, we
introduce a congestion price that measures
the marginal holding cost of additional
traffic. For model $m$ with a stable load
$\Lambda$, define
\begin{equation}
\Theta_m(\Lambda):=\Psi_m'(\Lambda)
=\frac{h\mu_m}{(\mu_m-\Lambda)^2},
\qquad 0\le\Lambda<\mu_m.
\label{eq:Theta}
\end{equation}
This congestion price can be decomposed as
\[
\Theta_m(\Lambda)
=\frac{h}{\mu_m-\Lambda}
+\frac{h\Lambda}{(\mu_m-\Lambda)^2}
\]
The first term represents the expected
holding cost of an additional visit at
the current load, including both waiting
and service. The second term captures the
congestion externality: the additional
holding cost imposed on other jobs as
their delays increase. Thus, the congestion
price accounts for the full system-wide
holding cost of sending additional traffic
to a model.

For each model $m$, let
$\theta_m:=\Theta_m(\Lambda_m)$ denote its
congestion price at the current load, and
define $\chi_m:=c_m+\theta_m$ as its
\emph{adjusted query cost}. This adjusted
cost combines the actual query fee $c_m$
with the marginal holding cost imposed on
the system. It therefore provides a measure
of the full marginal cost of querying model
$m$, which we use to guide routing decisions.
Under our assumptions, these costs are
model-specific but common across job types.
This approach follows the marginal-cost
pricing principle used in network routing
\citep{Kleinrock1964,BertsekasGallager1992}.

At a given vector of congestion prices
$\bm\theta$, the adjusted expected cost of
using tree $T$ for a type-$j$ job is
\begin{equation*}
A_{jT}(\bm\theta)
:=b_je_{jT}+\sum_{m=1}^M(c_m+\theta_m)v_{jTm}
=\omega_{jT}+\sum_{m=1}^M\theta_mv_{jTm}.
\end{equation*}
To see how it determines optimal routing,
consider a stable feasible mixture in which
tree $T$ is used with positive probability
$z_{jT}>0$. Suppose we shift a small probability
$\varepsilon>0$ from tree $T$ to another tree
$T'$ for the same job type. This changes the
routing assignments of type-$j$ jobs at rate
$\lambda_j\varepsilon$. For sufficiently small
$\varepsilon$, the new mixture remains feasible
and stable, and the first-order change in the
total cost rate is
$
\varepsilon\lambda_j
\bigl[A_{jT'}(\bm\theta)-A_{jT}(\bm\theta)\bigr].
$
If tree $T'$ has a lower adjusted cost than
tree $T$, this change is negative, so the
transfer reduces the total cost for sufficiently
small $\varepsilon$. Therefore, at an optimal
solution, every tree assigned positive
probability must minimize the adjusted cost
for its job type, evaluated at the congestion
prices induced by the optimal loads. Because
the master problem is convex by Theorem \ref{thm:convexity}, this condition
is also sufficient for global optimality.

\begin{theorem}[Congestion-price optimality conditions]
\label{thm:tree-KKT}
A stable feasible tree mixture $\bm z^*$ is optimal for \eqref{prob:master} if and only if, at its induced prices
$\theta_m^*=\Theta_m(\Lambda_m^*)$, $m\in[M]$, every positively used tree minimizes adjusted cost within its job type:
\begin{equation}
z_{jT}^*>0 \quad\Longrightarrow\quad A_{jT}(\bm\theta^*)
=\min_{T'\in\cT_L}A_{jT'}(\bm\theta^*).
\label{eq:tree-KKT}
\end{equation}
\end{theorem}

At fixed congestion prices, each job type can select its least-cost trees independently. To obtain a system optimum, however, the resulting mixture must generate the loads that produce those same prices through \eqref{eq:Theta}. Randomization among tied trees may be needed to meet this requirement. 

{\color{black}The optimality conditions identify which routing trees can be used, but do not construct them. The next section develops a Bayesian dynamic program that finds least-cost trees at given congestion prices without enumerating all possible trees. Its stopping and screening results further simplify this search.}

\section{Bayesian Policy Structure and Capacity Protection}
\label{sec:policy-structure}

We now turn to the routing and stopping
decisions of individual jobs. Recall from
Theorem~\ref{thm:tree-KKT} that, at the
optimal congestion prices, every tree used
by a job type must minimize its adjusted
expected cost. To identify these trees,
we develop a Bayesian dynamic program
that determines whether a job should stop
and make a classification or query another
model, based on the evidence collected so
far. We then characterize how the job's
confidence in its classification, the
informativeness of available models, and
query costs affect these decisions. These
results also allow us to derive explicit
bounds on the optimal traffic directed
to each model.

Throughout this section, we treat the adjusted query costs $\bm\chi=(c_m+\theta_m)_{m\in[M]} \in(0,\infty)^M$ as fixed. This allows us to analyze the routing and stopping decisions for each job type separately. We then apply these decisions to the master problem optimization in Section~\ref{sec:column-generation} to jointly determine the routing mixture and the congestion prices consistent with the resulting loads.

\subsection{The Bayesian state and the cost of stopping}

Under Assumption~\ref{ass:signals}, the information required for a job's remaining decisions can be summarized by $(p,S,r)$. Here, $p$ is its posterior probability of label $A$, $S$ is the set of unused models, and $r\in\{0,\ldots,|S|\}$ is the remaining query budget.

Although the choice of models may depend on
previously observed outputs, the selection
decisions themselves provide no additional
information about the true label once the
observed model-output pairs are known. This
follows because CAAR decisions depend only on
the observed history and independent
randomization. The following lemma formalizes
this property and characterizes both the
posterior probability of the true label and
the distribution of future model outputs.

\begin{lemma}[Bayesian updating under adaptive queries]
\label{lem:posterior-sufficiency}
At a reachable history $\mathfrak h_t=(m_s,x_s)_{s=1}^t$ of a type-$j$ job, the posterior probability of label $A$ is 
\[
p_j(\mathfrak h_t)
:=\frac{\rho_j\prod_{s=1}^t\ell_{jm_s}^A(x_s)}
{\rho_j\prod_{s=1}^t\ell_{jm_s}^A(x_s)+(1-\rho_j)\prod_{s=1}^t\ell_{jm_s}^B(x_s)}.
\]
Conditional on $Y_j=y$ and the observed history, unused model outputs remain mutually independent with likelihoods $\ell_{jm}^y$. Given the observed history, querying an unused model $m$ produces output $x$ with probability $s_{jm,x}(p_j(\mathfrak h_t))$ and updates the posterior to $\tau_{jm,x}(p_j(\mathfrak h_t))$. Conditioning additionally on the job's private random seed does not change these distributions.

\end{lemma}

The posterior determines the distribution of future outputs, while $S$ and $r$ determine which queries remain feasible. Two jobs with the same posterior can therefore face different decisions if their available models or budgets differ.
With repeated queries, every model remains available. Under the additional stationary conditional independence assumption in Appendix~\ref{sec:repeated-routing}, the state reduces to $(p,r)$.

We first evaluate the option of stopping at a state $(p,S,r)$. If the job is classified immediately, choosing label $A$ has expected error cost $b_j(1-p)$, while choosing label $B$ has expected error cost $b_jp$. The maximum-a-posteriori (MAP) decision selects the more likely label, with either label allowed at a tie. We therefore use MAP classification whenever a job stops. 

\begin{lemma}[MAP terminal decisions]
\label{lem:MAP-WLOG}
For any CAAR routing and stopping rule, replacing the terminal decisions by MAP decisions weakly lowers expected error cost and leaves query fees, model visits, queue loads, and holding costs unchanged.
\end{lemma}

The minimum expected cost of stopping is therefore $g_j(p):=b_j\min\{p,1-p\}$. To decide whether another query is worthwhile, we must compare this cost with the expected cost of querying and then making subsequent decisions optimally.

To evaluate the cost of future decisions at a
state $(p,S,r)$, consider the probability law
$\mathbb P_{j,p,S}$ under which the true label
is $A$ with probability $p$ and $B$ with
probability $1-p$. Conditional on the true
label, the outputs of the unused models in
$S$ retain their original likelihoods and
remain independent, as established in
Lemma~\ref{lem:posterior-sufficiency}. Let $\Pi^{c}_{j,S,r}$ %
denote the class of continuation policies that may stop
immediately or query at most $r$ distinct
models from $S$, with decisions based on
observed outputs and private randomization.
For any policy $\pi\in\Pi^{c}_{j,S,r}$,
let $N^\pi$ denote the number of future
queries, %
$M_t^\pi$ the model selected for
the $t$-th future query, and $\widehat Y^\pi$
the final classification. We define the
minimum expected remaining cost as
\begin{equation}
W_{j,S,r}(p;\bm\chi)
:=\min_{\pi\in\Pi^{c}_{j,S,r}}\mathbb E_{j,p,S}^{\pi}
\left[
\sum_{t=1}^{N^\pi}\chi_{M_t^\pi}
+b_j\mathbf1\{\widehat Y^\pi\ne Y_j\}
\right].
\label{eq:W-definition}
\end{equation}

This value accounts for the adjusted costs
of all future queries and the expected
penalty for the final classification.
Costs incurred by previous queries are
excluded because they cannot be changed
by future decisions. At fixed adjusted
query costs, randomization cannot improve
upon the best deterministic continuation
tree, since the expected cost of a randomized
policy is a weighted average of the costs
of deterministic trees. Thus, a deterministic
tree attains the minimum, and its terminal
decisions can be chosen according to the
MAP rule. After each query, the job faces
the same decision problem with an updated
posterior, a smaller set of available models,
and one fewer query remaining. This recursive
structure leads to the Bellman equation
developed next.

\subsection{Bellman decisions and consistency with system prices}\label{subsec:5.2}

We now use %
the remaining-cost function $W$ in \eqref{eq:W-definition} to formulate the Bellman equation. Given $(p,S,r)$, stopping has expected cost $g_j(p)$. On the other hand, %
querying model $m$ has adjusted cost $\chi_m$, followed by the cost of the decisions made after its output is observed. For $r\ge1$ and $m\in S$, define
\begin{equation*}
G_{jm,S,r}(p;\bm\chi)
:=\chi_m+\sum_{x\in\{A,B\}}s_{jm,x}(p)
W_{j,S\setminus\{m\},r-1}\bigl(\tau_{jm,x}(p);\bm\chi\bigr).
\end{equation*}
Here output $x$ occurs with probability $s_{jm,x}(p)$ and changes the posterior to $\tau_{jm,x}(p)$. Model $m$ is then removed from the available set, and the query budget falls by one. The function $G_{jm,S,r}$ therefore includes the cost of the next query and the optimal choices that follow it.

At arrival, a type-$j$ job has belief $\rho_j$, all models are available, and the query budget is $L$. Write $\Gamma_{j,L}(\bm\chi):=W_{j,[M],L}(\rho_j;\bm\chi)$ for its minimum adjusted expected cost. The following proposition gives the recursive comparison and connects this arrival value to the tree-selection problem.

\begin{proposition}[Bellman recursion and adjusted tree cost]
\label{prop:Bellman}
For every type $j$ and state $(p,S,r)$,
$W_{j,S,0}(p;\bm\chi)=g_j(p)$. For $r\ge1$,
\begin{equation}
W_{j,S,r}(p;\bm\chi)
=\min\left\{g_j(p),\ \min_{m\in S}G_{jm,S,r}(p;\bm\chi)\right\}.
\label{eq:Bellman}
\end{equation}
The value at arrival satisfies
\begin{equation}
\Gamma_{j,L}(\bm\chi)=W_{j,[M],L}(\rho_j;\bm\chi)
=\min_{T\in\cT_L}
\left\{b_je_{jT}+\sum_{m=1}^M\chi_mv_{jTm}\right\}.
\label{eq:root-tree-equivalence}
\end{equation}
\end{proposition}

The Bellman recursion can be solved by backward
induction, starting with $r=0$ and proceeding
to progressively larger query budgets. At each
state $(p,S,r)$, the job compares the cost of
stopping, $g_j(p)$, with the expected cost
$G_{jm,S,r}(p;\bm\chi)$ of querying each
available model. If stopping attains the
minimum, the job stops and makes a MAP
classification. Otherwise, it queries a model
with the lowest expected continuation cost
and repeats the same decision process after
observing its output. When multiple actions
attain the minimum, any of them may be chosen
without affecting the optimal expected cost.
We refer to these minimizing decisions,
including MAP classification upon stopping,
as \emph{Bellman-optimal actions}.

This recursion is a finite-horizon Bayesian
decision problem \citep{SmallwoodSondik1973},
in which choosing the next model corresponds
to selecting an experiment
\citep{Chernoff1959}. By
\eqref{eq:root-tree-equivalence}, the resulting
policy determines a least-cost routing tree
for each job type at the given congestion
prices. System optimality requires evaluating these decisions at the prices generated by the mixture itself. Combining Proposition~\ref{prop:Bellman} with Theorem~\ref{thm:tree-KKT} gives the following condition.

\begin{corollary}[Bellman implementation of a system optimum]
\label{cor:equilibrium-policy}
For a stable feasible mixture $\bm z$, set $\chi_m=c_m+\Theta_m(\Lambda_m(\bm z))$ for each $m\in[M]$.
The mixture is globally optimal if and only if, for every type $j$, each tree with $z_{jT}>0$ takes a Bellman-optimal action at every history it reaches with positive probability under these adjusted costs.
\end{corollary}

Corollary~\ref{cor:equilibrium-policy} connects the job-level Bellman decisions to system-wide optimality. The recursion identifies least-cost trees at given prices, but their mixture is optimal only when its induced loads generate those same prices. Several trees
may minimize adjusted cost for a given job type while generating
different model loads. Randomization among these tied trees may
therefore be necessary to obtain a mixture consistent with the
prices used to select them. The Bellman-optimality requirement
applies only to histories that a positively used tree reaches with
positive probability. Actions at unreachable histories affect
neither its expected cost nor the traffic it generates.

The Bayesian dynamic program is also
computationally manageable for modest numbers
of models. Under conditional independence,
histories that contain the same model-output
pairs, even in different orders, lead to the
same posterior and continuation problem.
Consequently, the number of states that
must be considered is substantially smaller
than the number of possible query histories.
For a fixed prior and depth limit $L$, the
number of states per job type is at most
$\sum_{d=0}^{L}\binom{M}{d}2^d\le3^M$,
including terminal states. For example,
with ten models, there are at most 1,161
states when $L=3$, and 59,049 even when
all ten models may be queried. Moreover,
at fixed congestion prices, the dynamic
programs can be solved independently
across job types, and states within each
backward-induction layer can be evaluated
in parallel. The stopping rules and model
screening conditions developed next can
further reduce computation by identifying
states where no additional query is needed
and models that need not be considered.

\subsection{Two-sided stopping and the effect of query costs}

We now examine how the optimal stopping decision depends on the posterior belief and adjusted query costs. We show that stopping is optimal when the posterior strongly favors either label and that increasing adjusted query costs weakly expands the stopping region.

The key is to characterize how the remaining-cost function varies with the posterior belief. Fix $j,S,r,\bm\chi$ and consider a deterministic continuation policy $\pi\in\Pi^{c}_{j,S,r}$ with fixed routing, stopping, and terminal classification rules. For $y\in\{A,B\}$, let $W^{\pi,y}(\bm\chi)$ denote its expected remaining cost conditional on the true label being $y$. Because the policy and the conditional signal distributions are fixed, these conditional costs do not depend on the posterior belief $p$. Consequently, the expected remaining cost under policy $\pi$ is

\begin{equation*}
W^\pi(p;\bm\chi)
=pW^{\pi,A}(\bm\chi)+(1-p)W^{\pi,B}(\bm\chi),
\end{equation*}
which is affine in $p$. Since there are finitely many deterministic continuation trees and randomization cannot improve upon their minimum, the optimal remaining-cost function $W_{j,S,r}(p;\bm\chi)$ is the pointwise minimum of finitely many affine functions. This observation yields the following properties.

\begin{lemma}[Property of the remaining cost function]
\label{lem:posterior-geometry}
For fixed $j,S,r,\bm\chi$, the function $W_{j,S,r}(\cdot;\bm\chi)$ is continuous, concave, and piecewise affine with finitely many pieces in $p$. Moreover,
$0\le W_{j,S,r}(p;\bm\chi)\le g_j(p)$ and $ W_{j,S,r}(0;\bm\chi)=W_{j,S,r}(1;\bm\chi)=0$.
For $r\ge1$ and $m\in S$, $G_{jm,S,r}(\cdot;\bm\chi)$ is also continuous, concave, and piecewise affine with finitely many pieces.
\end{lemma}

Since immediate stopping is always feasible, $W_{j,S,r}(p;\bm\chi)\le g_j(p)$. Stopping is optimal exactly when equality holds. We therefore define the stopping region as
\[
\mathcal Z_{j,S,r}(\bm\chi)
:=\{p\in[0,1]:W_{j,S,r}(p;\bm\chi)=g_j(p)\}.
\]
By Lemma~\ref{lem:posterior-geometry}, $W_{j,S,r}(p;\bm\chi)$ is concave in $p$, whereas $g_j(p)$ is affine on each of the intervals $[0,1/2]$ and $[1/2,1]$. Consequently, the difference $g_j(p)-W_{j,S,r}(p;\bm\chi)$ is nonnegative and convex on each interval. Moreover, this difference vanishes at both $p=0$ and $p=1$. Its zero set on $[0,1/2]$ must therefore be an interval starting at $0$, while its zero set on $[1/2,1]$ must be an interval ending at $1$. This establishes the two-sided threshold structure of the optimal stopping region.

\begin{proposition}[Two-sided posterior stopping]
\label{prop:stopping-threshold}
For every $j,S,r,\bm\chi$, there are cutoffs $0\le p^-_{j,S,r}(\bm\chi)\le1/2\le p^+_{j,S,r}(\bm\chi)\le1$ such that
$\mathcal Z_{j,S,r}(\bm\chi)
=[0,p^-_{j,S,r}(\bm\chi)]\cup[p^+_{j,S,r}(\bm\chi),1]$.
If $r\ge1$, define $\underline\chi_S:=\min_{m\in S}\chi_m$. The cutoffs satisfy
\begin{equation*}
p^-_{j,S,r}(\bm\chi)\ge\min\left\{\frac12,\frac{\underline\chi_S}{b_j}\right\},
\qquad
p^+_{j,S,r}(\bm\chi)\le1-\min\left\{\frac12,\frac{\underline\chi_S}{b_j}\right\}.
\end{equation*}
\end{proposition}

Proposition~\ref{prop:stopping-threshold} establishes a two-sided threshold structure for optimal stopping. When the posterior is below the lower cutoff, the job stops and selects label $B$; when it is above the upper cutoff, the job stops and selects label $A$. Between the cutoffs, querying another model is strictly preferable to stopping. If this interval is nonempty, stopping and querying are tied at either cutoff. If both cutoffs equal $1/2$, stopping is optimal for every posterior belief, but need not tie with querying. This structure is consistent with the classical two-sided stopping rules in Bayesian sequential testing \citep{Wald1947,ArrowBlackwellGirshick1949,DeGroot1970}.

The proposition also provides an explicit sufficient condition for stopping without solving the Bellman recursion. Since any additional query incurs an adjusted cost of at least $\underline\chi_S$, stopping is optimal whenever $g_j(p)=b_j\min\{p,1-p\}\le\underline\chi_S$. Thus, when the posterior sufficiently favors either label, further information acquisition cannot justify even the cost of the next query. Outside this guaranteed stopping region, querying may still be suboptimal because the expected information benefit need not outweigh its adjusted cost.

We next examine how optimal stopping decisions respond to changes in adjusted query costs, holding the job type, signal distributions, available models, and remaining query budget fixed. 

\begin{corollary}[Higher query costs expand stopping]
\label{cor:cost-stopping}
If $\bm\chi'\ge\bm\chi$ componentwise on $S$, then, for every $p\in[0,1]$,
$W_{j,S,r}(p;\bm\chi')\ge W_{j,S,r}(p;\bm\chi)$ and $
\mathcal Z_{j,S,r}(\bm\chi)\subseteq\mathcal Z_{j,S,r}(\bm\chi')$.
\end{corollary}

Increasing query costs weakly raises the expected cost of every fixed continuation policy, while leaving the cost of immediate stopping unchanged. 
Higher query costs thus preserve every previously optimal stopping decision, and hence, increasing one model's congestion price while holding the other prices fixed can only expand the stopping region. Changes in arrival rates or service capacities may alter congestion prices and the histories reached, so the corollary does not by itself describe the full system response.

\subsection{Information-value screening and capacity protection}
\label{subsec:information-cost-frontier}

Upon choosing the next model to query, the Bellman recursion requires comparing all available models at each decision state. We now derive an information-value bound that allows us to exclude models whose potential benefit cannot justify their adjusted query costs. This bound first provides a screening rule for selecting the next model at a given posterior belief. We then use it to identify models that can be excluded from the optimization altogether and to derive explicit upper bounds on their optimal loads.

The benefit of an additional query arises from the information it provides, which changes the posterior belief and potentially reduces subsequent error probability. To bound this benefit, we first establish how sensitive the optimal remaining cost is to changes in the posterior. The following lemma provides a Lipschitz bound, using the concavity and value bounds established in Lemma~\ref{lem:posterior-geometry}.

\begin{lemma}[Sensitivity of remaining cost to belief]
\label{lem:W-lipschitz}
For every $j,S,r,\bm\chi$ and all $p,p'\in[0,1]$,
$\left|W_{j,S,r}(p;\bm\chi)-W_{j,S,r}(p';\bm\chi)\right|\le b_j|p-p'|$.
\end{lemma}

To relate the Lipschitz bound to signal quality, we assume that output $A$ is more likely under label $A$ than under label $B$. Note that the preceding Bayesian and stopping results do not require this condition.

\begin{assumption}[Positive signal orientation]
\label{ass:orientation}
For every $j\in[J]$ and $m\in[M]$, $q_{jm}^A+q_{jm}^B>1$.
\end{assumption}

Define $d_{jm}:=q_{jm}^A+q_{jm}^B-1$. Under Assumption~\ref{ass:orientation}, $d_{jm}\in(0,1)$ equals the total variation distance between the conditional output distributions:
$
d_{jm}=\frac12\sum_{x\in\{A,B\}}
|\ell_{jm}^A(x)-\ell_{jm}^B(x)|.
$
Thus, smaller $d_{jm}$ indicates a less informative model. Let $P'$ denote the posterior after querying model $m$ at belief $p$. By \eqref{eq:predictive}--\eqref{eq:update},
$\mathbb E[|P'-p|]
=\sum_{x\in\{A,B\}}s_{jm,x}(p)
\left|\tau_{jm,x}(p)-p\right|
=2p(1-p)d_{jm}.$
Each output contributes $p(1-p)d_{jm}$ to this expectation, even with asymmetric accuracies. Hence the expected posterior change is small when the signal is weak or the current belief is close to $0$ or $1$. Together with Lemma~\ref{lem:W-lipschitz}, this identity bounds the potential benefit of an additional query.

To bound the benefit of querying model $m$, consider the feasible alternative of skipping $m$ and using at most $r-1$ queries from $S\setminus\{m\}$. Both choices leave the same models and query budget for subsequent decisions, but querying $m$ provides an additional signal. Lemma~\ref{lem:W-lipschitz} and the posterior-change formula bound the resulting reduction in expected future cost by $2b_jp(1-p)d_{jm}$. Moreover, querying $m$ cannot be optimal if its adjusted cost exceeds the immediate stopping cost $g_j(p)$. Combining these bounds yields the following necessary condition.

\begin{theorem}[Statewise information-cost frontier]
\label{thm:statewise-information-frontier}
Under Assumption~\ref{ass:orientation}, if querying $m\in S$ is Bellman optimal at $(p,S,r)$ with $r\ge1$, then
\begin{equation}
\chi_m\le b_j\min\{p,1-p,2p(1-p)d_{jm}\}.
\label{eq:statewise-frontier}
\end{equation}
\end{theorem}

\begin{figure}[t]
\centering
\includegraphics[width=0.8\textwidth]{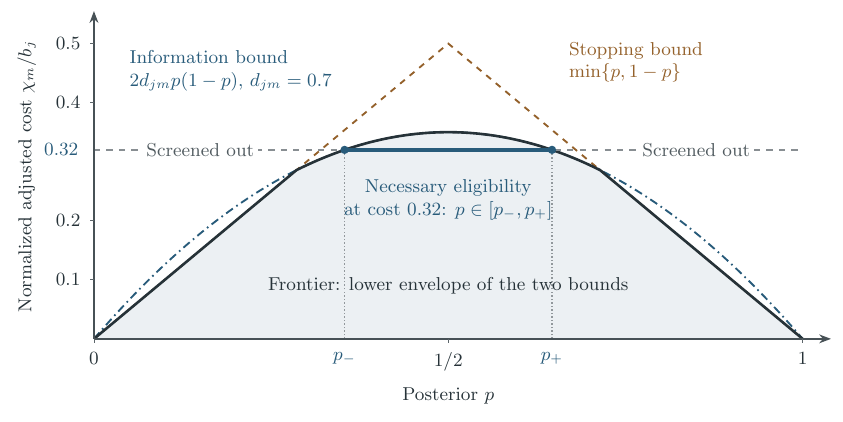}
\caption{The next-query screening bound for $d_{jm}=0.7$. The solid boundary is the lower envelope of the stopping bound (ochre, dashed) and the information bound (blue, dash-dot). Shading marks cost-posterior pairs satisfying the necessary condition \eqref{eq:statewise-frontier}. At normalized adjusted cost $\chi_m/b_j=0.32$, the model remains a candidate only for $p\in[p_-,p_+]\approx[0.354,0.646]$, including the endpoints. Passing the screen does not establish optimality of querying.}
\label{fig:information-frontier}
\end{figure}

A model violating this condition cannot be an optimal next query, even in a tie. Define the remaining candidate set by
$\cF_{j,S}(p;\bm\chi):=\{m\in S:\chi_m\le b_j\min\{p,1-p,2p(1-p)d_{jm}\}\}.$
Restricting the next-query comparison to this set leaves the Bellman value unchanged:
\begin{equation}
W_{j,S,r}(p;\bm\chi)=\min\left\{g_j(p),\ \min_{m\in\cF_{j,S}(p;\bm\chi)}G_{jm,S,r}(p;\bm\chi)\right\},\qquad r\ge1.
\label{eq:Bellman-screened}
\end{equation}
The inner minimum is $+\infty$ when the candidate set is empty, so stopping is then optimal. This restriction applies only to the next action. An excluded model remains available and may become eligible after another output changes the posterior. A model in this set is only a candidate, and the recursion still compares it with stopping and the other candidates. Figure~\ref{fig:information-frontier} illustrates how the screen combines the stopping-cost and information-value bounds to determine eligible models.

The same bound also lets us rule out some models before solving for a policy or its congestion prices. Since $2p(1-p)\le1/2$, Theorem~\ref{thm:statewise-information-frontier} limits the adjusted cost of any optimal type-$j$ query to $b_jd_{jm}/2$. Meanwhile, model $m$'s adjusted cost is at least its zero-load value, $\chi_m^{\mathrm{idle}}:=c_m+\Theta_m(0)=c_m+h/\mu_m$. Even an idle model incurs holding cost during service. Positive use creates congestion and raises its adjusted cost strictly above this minimum. Define
\begin{equation}
\cF_j^0:=\left\{m\in[M]:\chi_m^{\mathrm{idle}}<\frac{b_jd_{jm}}{2}\right\},\qquad \cF^0:=\bigcup_{j=1}^J\cF_j^0.
\label{eq:exact_screen_f_0}
\end{equation}
These sets can be computed from the upfront given problem parameters. Also let $\overline v_m:=\max_{j\in[J]}b_jd_{jm}/2$ denote the largest of these information-value bounds across job types. 
\begin{corollary}[Model exclusion and optimal load bounds]
\label{cor:exogenous-frontier}
Under Assumption~\ref{ass:orientation}, every optimal CAAR tree-mixture $\bm z^*$ satisfies:
\begin{enumerate}[label=(\roman*),leftmargin=*]
\item If $m\notin\cF_j^0$, type $j$ queries $m$ with zero probability. Restricting each type's tree space to $\cF_j^0$ preserves the optimal value.
\item If $m\notin\cF^0$, then $\Lambda_m^*=0$.
\item If $\overline v_m\le\chi_m^{\mathrm{idle}}$, then $\Lambda_m^*=0$. If $\Lambda_m^*>0$, then $\overline v_m>\chi_m^{\mathrm{idle}}$ and
\begin{equation}
\Lambda_m^*\le\mu_m-\sqrt{\frac{h\mu_m}{\overline v_m-c_m}}.
\label{eq:information-load-cap}
\end{equation}
\end{enumerate}
\end{corollary}

The strict inequality in \eqref{eq:exact_screen_f_0} distinguishes this test from the screen at a given belief. If a model's idle cost already equals its information-value cap, any positive traffic makes its adjusted cost exceed that cap. Such a model cannot receive positive traffic from that type at the system optimum. A model that passes the test may still remain unused because of competing models or congestion.

Part~(iii) expresses the same restriction as a load bound. For an active model, the adjusted cost cannot exceed $\overline v_m$, whereas its marginal holding cost grows without bound as load approaches capacity. Hence its optimal load leaves the explicit capacity buffer
\begin{equation*}
\mu_m-\Lambda_m^*\ge\sqrt{\frac{h\mu_m}{\overline v_m-c_m}}>0.
\end{equation*}
\textcolor{black}{
This capacity buffer highlights why maximizing server utilization need not minimize total cost. For an active model, the lower bound on spare capacity increases with the holding cost $h$ and query fee $c_m$, but decreases with the information-value bound $\overline v_m$, holding other parameters fixed. Because arrival rates do not enter this bound, higher demand cannot eliminate the requirement to maintain spare capacity. These comparisons concern the bound itself, rather than the actual optimal loads.}

The two screening rules serve complementary purposes: the statewise rule excludes suboptimal next queries at a given belief and price, while the exogenous rule removes models from the policy search altogether. \textcolor{black}{Appendix~\ref{apdx:policy-special-cases} presents a water-filling solution for the one-query ($L=1$) case and a complete policy characterization for the two-model case. In the next section, we use the Bayesian recursion as a pricing subroutine within column generation to jointly compute the optimal routing mixture and its congestion prices.}

\FloatBarrier
\section{Column Generation for Computing an Optimal CAAR Policy}
\label{sec:column-generation}
The preceding section identifies least-cost routing trees at given congestion prices. We now jointly compute the optimal routing mixture and the congestion prices induced by its model loads as solving \eqref{prob:master} by the column generation and the Bayesian recursion. Although the master problem is convex, the number of feasible trees can make direct enumeration impractical. Column generation treats each routing tree as a column and optimizes a mixture over a growing collection of trees. At the resulting congestion prices, the Bayesian recursion searches for improving trees to add. An explicit optimality-gap bound provides a stopping criterion before exact optimality is reached. Assumptions~\ref{ass:signals} and \ref{ass:orientation} remain in force.

Fix a depth limit $K\in\{0,\ldots,M\}$. We first use the model exclusion result in Corollary~\ref{cor:exogenous-frontier}, and restrict the possible tree mixture to the collection of $\cF_j^0$ for each job type $j$. Let $\cT_{j,K}$ contain the trees with at most $K$ queries using models in $\cF_j^0$ and denote the system optimal cost $C^*_K$.

The master problem determines how each job type's traffic should be divided among routing trees. Solving it directly would require considering all feasible trees at once. The \emph{restricted master} makes the same allocation decision using only the trees already generated for each type. Denote this set by $\cC_j\subseteq\cT_{j,K}$. Initially, $\cC_j$ contains only the tree that stops at the root and chooses the more likely label under the prior. Given these alternatives, the restricted master chooses their routing mixture:
\begin{equation}
\label{prob:cg-master}
\begin{aligned}
U=\min_{\bm z}\quad&
\sum_j\lambda_j\sum_{T\in\cC_j}z_{jT}\omega_{jT}
+\sum_m\Psi_m(\Lambda_m(\bm z))\\
\text{s.t.}\quad&
\sum_{T\in\cC_j}z_{jT}=1,\qquad z_{jT}\ge0,
\qquad j\in[J],\ T\in\cC_j.
\end{aligned}
\end{equation}
The loads follow \eqref{eq:load-tree}, with zero weight on trees outside $\cC_j$. This convex problem accounts for the congestion created by the combined traffic assigned to the available trees. The stopping trees ensure a finite-cost feasible solution. Its optimal value $U$ is an upper bound on $C_K^*$.

The optimal loads in \eqref{prob:cg-master} determine congestion prices $\theta_m=\Theta_m(\Lambda_m)$ and adjusted query costs $\chi_m=c_m+\theta_m$. At these prices, define
\begin{equation}
\nu_j:=\min_{T\in\cC_j}A_{jT}(\bm\theta),
\qquad
\phi_j(\bm\theta):=\min_{T\in\cT_{j,K}}A_{jT}(\bm\theta).
\label{eq:cg-pricing}
\end{equation}
Every positively used tree in the restricted-master solution has adjusted cost $\nu_j$. If $\phi_j(\bm\theta)<\nu_j$, shifting a small amount of traffic to a minimizing tree reduces total cost. The search for this tree is the \emph{tree-pricing problem}. It takes congestion prices as inputs, rather than optimizing those prices directly. Proposition~\ref{prop:Bellman} solves the search in \eqref{eq:cg-pricing} by Bayesian recursion:
\[
\phi_j(\bm\theta)
=W_{j,\cF_j^0,K_j}(\rho_j;\bm\chi),
\qquad K_j:=\min\{K,|\cF_j^0|\}.
\]
Specifically, when solving the recursion, we can use the state-dependent screening results in Proposition~\ref{prop:stopping-threshold} and Theorem~\ref{thm:statewise-information-frontier} to further speed up the process of finding the next model to query or stop. 
The recursion also returns a minimizing tree that can be included in the next iteration. 

Starting with the immediate-stopping trees, column generation repeats:
\begin{enumerate}[label=\textbf{Step \arabic*.},leftmargin=*]
\item Solve the restricted master and compute its induced prices and the values $\nu_j$.
\item For each type $j$, use the Bayesian recursion to obtain a minimizing tree $T_j$ and its adjusted cost $\phi_j(\bm\theta)$.
\item If $\phi_j(\bm\theta)=\nu_j$ for every type, return the current mixture. Otherwise, update $\cC_j\leftarrow\cC_j\cup\{T_j\}$ for each type with $\phi_j(\bm\theta)<\nu_j$ and repeat.
\end{enumerate}
Tree pricing identifies a routing alternative, but not how much traffic it should receive. The next restricted master chooses its weight jointly with those of the existing trees and updates congestion prices. Retaining all generated trees ensures that the restricted-master value cannot increase.

The same tree search measures how much the current policy can still improve. For type $j$, $\nu_j-\phi_j(\bm\theta)$ is the largest reduction in adjusted tree cost at the current prices. Replacing each convex holding cost by its tangent at the current load gives a lower bound on system cost. Weighting the adjusted-cost differences by arrival rates therefore yields the following guarantee.

\begin{proposition}[Optimality gap and finite convergence]
\label{thm:cg-certificate}
\label{thm:cg-correctness}
At an exact restricted-master optimum with value $U$ and induced prices $\bm\theta$, define
\begin{equation}
LB:=U-\sum_j\lambda_j\bigl(\nu_j-\phi_j(\bm\theta)\bigr).
\label{eq:cg-gap-bound}
\end{equation}
Then $LB\le C_K^*\le U$. The restricted-master solution is globally optimal if and only if $\phi_j(\bm\theta)=\nu_j$ for every type. For fixed $K$, column generation with exact master and tree-pricing solves terminates after finitely many column additions at $C_K^*$.
\end{proposition}

With exact master and tree-pricing solutions and all generated trees retained, column generation reaches $C_K^*$ from any initial tree set admitting a stable feasible mixture. The initial set can affect the number of iterations, but not the optimal value, and the immediate-stopping trees guarantee a feasible start. Finite termination does not imply a polynomial running time, since the number of possible trees can be large. The lower bound also supports early stopping. For a cost tolerance $\varepsilon\ge0$, stopping when $U-LB\le\varepsilon$ guarantees a policy within $\varepsilon$ of the depth-$K$ optimum. The bound in \eqref{eq:cg-gap-bound} accounts for all remaining improvements, including reallocating traffic and the resulting changes in congestion prices. The next section bounds the separate cost of limiting query depth. This can justify a smaller column-generation problem within a prescribed tolerance for the additional cost.

\FloatBarrier
\section{The Value of Depth and Effective Depth}
\label{sec:effective-depth}

The preceding algorithm optimizes routing at a chosen query depth, and we now investigate how much can be lost by choosing a smaller depth limit. Additional queries can improve classification, but also incur fees and increase congestion. In this section, we derive an upper bound on the cost difference given a depth limit. This bound gives a sufficient depth for a prescribed cost tolerance using the model parameters before optimization.

To begin with, a query is more informative when its output distributions under the two labels have less overlap. For type $j$ and model $m$, we measure this overlap by $\beta_{jm}$, the \emph{Bhattacharyya coefficient} \citep{bhattacharyya1943measure}, and express its information contribution as $\mathcal I_{jm}$:
\begin{equation*}
\beta_{jm}:=\sum_{x\in\{A,B\}}
\sqrt{\ell_{jm}^A(x)\ell_{jm}^B(x)},
\qquad
\mathcal I_{jm}:=-\log\beta_{jm}.
\end{equation*}
Under Assumptions~\ref{ass:signals} and \ref{ass:orientation}, we have $0<\beta_{jm}<1$. {\color{black}This coefficient measures how similar the model's output distributions are under the two labels. Smaller overlap means that the output distinguishes the labels more clearly, giving a larger information contribution $\mathcal I_{jm}$. Taking negative logarithms converts products of overlap coefficients into sums of information contributions.}

To cover every adaptive route, order the contributions as $\mathcal I_{j,(1)}\le\cdots\le\mathcal I_{j,(M)}$, and let $\mathcal I_j(r):=\sum_{k=1}^r\mathcal I_{j,(k)}$, with $\mathcal I_j(0):=0$. This sum uses the $r$ least informative models, and any other set of $r$ models has at least as large a sum. Equivalently, $e^{-\mathcal I_j(r)}$ is the product of the $r$ largest overlap coefficients. Note that this ordering is used to evaluate the bound, but not a query order.

The potential error from stopping after $K$ queries leads to the following bound on the increase in system cost. Let $B_j:=b_j\sqrt{\rho_j(1-\rho_j)}$ and define
\begin{equation*}
\mathcal E(K):=
\begin{cases}
\displaystyle\sum_{j=1}^J\lambda_jB_j e^{-\mathcal I_j(K)}, &0\le K<M,\\
0,&K=M.
\end{cases}
\end{equation*}
The weights account for arrival rates, error penalties, and prior uncertainty. The endpoint $\mathcal E(M)=0$ records that no depth truncation occurs when all $M$ queries remain available. All quantities in the bound come from the given model and job parameters, so $\mathcal E(K)$ can be evaluated before optimization. 

\begin{theorem}[Cost of limiting adaptive depth]
\label{thm:adaptive-depth-truncation}
For $0\le K<L\le M$,
$0\le C_K^*-C_L^*\le\mathcal E(K)$.
\end{theorem}

The comparison in Theorem \ref{thm:adaptive-depth-truncation} allows the routing mixture and its induced congestion prices to be reoptimized at each depth. To construct the bound, take an optimal depth-$L$ tree mixture and stop each branch after at most $K$ queries, using MAP wherever the original tree would continue. This removes visits, so every model load and both query and holding costs weakly decrease, but the classification loss can increase. The truncated mixture is a feasible depth-$K$ policy, which is sufficient to bound the optimal depth-$K$ cost.

In particular, the value of one more layer satisfies $0\le C_L^*-C_{L+1}^*\le\mathcal E(L)$ for $L<M$. The bound is conservative because it uses the least informative models and discards the fee and congestion savings from truncation. It measures the possible cost benefit of allowing deeper branches, rather than how often those branches will be used. %

The truncation bound then gives a sufficient depth for a prescribed tolerance. For a cost tolerance $\varepsilon>0$ per unit time, choose the smallest depth for which this bound is small enough: 
$
K_\varepsilon:=
\min\{K\in\{0,\ldots,M\}:\mathcal E(K)\le\varepsilon\}.
$
This set is nonempty because $\mathcal E(M)=0$. The choice $K=0$ means classifying every job from its prior without a query, whereas $K=M$ recovers the full-depth CAAR problem. The chosen depth is a cap, and individual jobs can still stop earlier.

\begin{corollary}[A sufficient depth for a cost tolerance]
\label{cor:adaptive-epsilon-depth}
Any optimal depth-$K_\varepsilon$ policy satisfies
$0\le C_{K_\varepsilon}^*-C_M^*\le\varepsilon$.
\end{corollary}

The choice $K_\varepsilon$ uses the full collection of information contributions. A simpler sufficient depth uses only the smallest contribution $\underline{\mathcal I}:=\min_{j,m}\mathcal I_{jm}>0$ and the total weight $B_\Sigma:=\sum_j\lambda_jB_j$. Set
\[
\widehat K_\varepsilon:=
\min\left\{M,\max\left\{0,
\left\lceil\frac{\log(B_\Sigma/\varepsilon)}
{\underline{\mathcal I}}\right\rceil\right\}\right\}.
\]
This depth limit also satisfies the guarantee in Corollary~\ref{cor:adaptive-epsilon-depth}, because $\mathcal I_j(K)\ge K\underline{\mathcal I}$ gives $\mathcal E(K)\le B_\Sigma e^{-K\underline{\mathcal I}}$ for $K<M$, and the gap is zero at $K=M$. Stronger signals permit a smaller sufficient depth, while higher arrival rates or penalties require tighter control of per-job error to meet the same cost-rate tolerance. With the other parameters fixed, the sufficient depth grows at most logarithmically in $1/\varepsilon$, up to the cap $M$.

\section{Comparison with Benchmark Routing Policies}
\label{sec:numerical}
This section separates the value of output-dependent model choice from the value of accounting for congestion. The comparisons are motivated by fixed-order cascades \citep{chen2023frugalgpt,jiang2026cascadia} and adaptive cascade routing with exogenous query costs \citep{dekoninck2025unified}. Fixed-ranking (FR) and mixed-ranking (MR) policies restrict model order but optimize the same congestion-aware objective as CAAR. Fee-only (FO) and service-time-only (SO) policies retain output-adaptive model choice but omit load-dependent waiting costs. %

\subsection{Fixed rankings and the value of adaptive model choice}

Under FR, each job type follows one optimized model order. In many router applications, models are ranked by their ability or quality, and each job are routed in a fixed order. Under MR, each arrival draws an order from an optimized type-specific distribution before observing outputs. Both FR and MR allow randomized, history-dependent stopping, including immediate stopping. A continuing job must query the next model in its assigned order. All three classes minimize the same total cost, so $C_K^*\le C_{\mathrm{MR}}^*\le C_{\mathrm{FR}}^*$. MR can distribute traffic across orders to reduce congestion. CAAR additionally chooses the next model after seeing the preceding outputs. For example, a screener's output can determine which of two specialists to consult. We can construct such a family with three models that the cost ratio between the fixed-cascade benchmarks and CAAR goes to infinity. 

\begin{proposition}[The cost of a predetermined model order]
\label{prop:fixed-ranking-gap}
For any fixed $h>0$, there exists a family of one-type instances with three models and depth $K=2$ such that
$
\frac{C_{\mathrm{FR}}^*}{C_2^*}
\ge \frac{C_{\mathrm{MR}}^*}{C_2^*}
\to \infty.
$
The arrival rate, positive query fees, and service rates are fixed across the family. Immediate stopping is allowed, and every admissible policy is stable in every instance.
\end{proposition}

In the construction, signals become more accurate while the error penalty increases relative to fixed query fees. CAAR uses the screener's output to select the more informative specialist. A predetermined pair cannot reproduce this choice. The divergence therefore changes both signals and penalties, rather than increasing the penalty alone. The two-query cap is important for this construction. Longer query chains can consult both specialists and stop when the evidence is sufficient, while randomizing across rankings can further distribute traffic.

For a finite-penalty comparison, consider a single job type ($J = 1$) and four models ($M = 4$) with $L=2$, prior $0.5$, $b=2$, and $h=0.01$. Base mean accuracies are $(0.70,0.70,0.80,0.80)$, service rates are $1.25$, and fees are $0.01$. Five independent parameter draws perturb each mean accuracy uniformly by up to $0.01$, each service rate by up to $6\%$, and each fee by up to $10\%$. For perturbed means $a_m$, set $q_m^A=a_m+\delta s_m$ and $q_m^B=a_m-\delta s_m$, where $s=(0,0,1,-1)$ and $\delta$ represents the level of model specialization. The last two models specialize in opposite labels, without changing their balanced-prior single-query accuracy. We use specializations $\delta\in\{0,0.025,0.05,0.075,0.10\}$ and offered loads $\alpha=\lambda/\min_m\mu_m\in\{0.5,0.7,0.9\}$. The whole setup gives a total of 75 instances. Each draw stays fixed across this grid.

The CAAR optimal policy is solved with the Bayesian recursion and Section~\ref{sec:column-generation}'s algorithm. FR optimizes stopping for each of the 12 ordered model pairs, and MR optimizes their mixtures. Costs are evaluated analytically. Let $\pi_B$ denote an optimal policy computed under benchmark $B$'s
objective and policy restrictions. We write $C_B^*:=C_L(\pi_B)$ for
its actual total cost, including holding costs. Throughout this section, CAAR savings against benchmark $B$ is computed by $(C_B^*-C_L^*)/C_B^*\times 100\%$. All policies allow immediate stopping. All instances satisfy $\sum_j\lambda_j<\min_m\mu_m$, so every admissible policy is stable.

\begin{table}[htbp]
\small
\begin{tabular}{llrrrr}
\toprule
Comparison & Benchmark & Types & Query cap & Mean & Maximum\\
\midrule
Fixed cascade order & FR & 1 & 2 & 8.07 & 20.57\\
            & MR & 1 & 2 & 5.08 & 14.43\\
Congestion awareness & FO & 1 & 3 & 7.55 & 24.63\\
              & SO & 1 & 3 & 7.30 & 24.63\\
              & FO & 2 & 3 & 7.06 & 24.71\\
              & SO & 2 & 3 & 5.80 & 24.42\\
\bottomrule
\end{tabular}
\centering
\caption{CAAR savings relative to each benchmark in the specialist setup (\%).} %
\label{tab:policy-comparisons}
\end{table}

Table~\ref{tab:policy-comparisons} shows noticeable gains at a fixed penalty. At the strongest specialization, average savings are 15.97\% against FR and 13.23\% against MR. Higher load increases FR's disadvantage in this grid. Its mean gap rises from 5.34\% at $\alpha=0.5$ to 12.33\% at $\alpha=0.9$, with a maximum of 20.57\%. MR already distributes traffic, and its mean gap stays near 5\% across these loads. %

\subsection{The value of congestion awareness}

We now retain the same adaptive policy class and change the objective used to design it. FO minimizes classification penalties and query fees, $S_L(\pi)=\sum_j\lambda_j[b_je_j(\pi)+\sum_m c_mv_{jm}(\pi)]$. SO adds the expected processing cost, and minimizes  $S_L(\pi) + h\sum_m\Lambda_m(\pi)/\mu_m$. Thus SO charges $c_m+h/\mu_m$ per query, independently of the traffic that the policy generates. It can be seen as a simple exogenous latency cost. CAAR also accounts for waiting, whose cost depends on all types' routing decisions. All three policies consider adaptive routing with no fixed model orders.

We consider three queries, $L=3$, in the cases examined below containing three different scenarios for $M=4$. The main scenario retains the \emph{specialists setup} as the preceding parameters and grid for fixed cascade, but considers both one and two job types. In the two-type version, each of the two job types accounts for half of the load and their mean accuracies are $a_{1m}=a_m+0.025s_m$ and $a_{2m}=a_m-0.025s_m$. This gives a total of 150 instances for the main scenario, that is, 5 draws $\times 5$ specialization levels $\times 3$ loads $\times 2$ type configurations. %

We further consider two additional setup designs to separate the mechanisms behind. A \emph{shared-server design} has two job types, $b_1=b_2=1$, $h=0.03$, service rates $(1,1.25,1.25,1.5)$, and fees $(0.01,0.012,0.012,0.008)$. Both output accuracies are $(0.84,0.80,0.70,0.72)$ for type 1 and $(0.84,0.70$, $0.80,0.72)$ for type 2. We fix $\lambda_2=0.15$ and vary $\lambda_1$ over $0.05,0.25,0.45,0.65,0.80$, for the base case and five draws with the preceding perturbation ranges, keeping $\mu_1=1$. Each model receives a separate accuracy perturbation, applied equally to both job types. 
A \emph{symmetric design} uses four identical models with accuracy $0.70$, service rate 1, fee $0.005$, and $h=0.03$. We use total arrival rates $0.1,0.3,0.5,0.7,0.9$. The one-type cases have penalty $b=0.5$ or 1. In the two-type cases, each type receives half the total arrivals, with penalties $b_1= 0.5$, $b_2=1$. All priors remain $0.5$. For brevity, we exclude the results for these two additional setups from Table 1 and discuss selected instances that illustrate the key findings.

In computation, we optimize mixtures over all 3,137 query/stop trees per type and check their linear cost minima with dynamic programming. FO and SO first minimize their own objectives. Among all tied minimizers, we choose the mixture with the lowest actual total cost, jointly across types. This favors the benchmarks and avoids attributing a gap to poor tie-breaking. All reported solutions are within $10^{-4}C_{\mathrm{stop}}$ of the relevant optimum, where $C_{\mathrm{stop}}=\sum_j\lambda_jb_j\min\{\rho_j,1-\rho_j\}$.

The specialist results reported in Table~\ref{tab:policy-comparisons} show that congestion awareness remains valuable with three queries and two types. Considering service-time fees reduces the mean gap compared to FO, but it does not always improve actual cost. This is because a slightly faster model can attract more traffic and create enough waiting to outweigh the processing-time saving. Both benchmarks can therefore benefit a lot from explicitly pricing the loads they induce.

The shared-server cases show how traffic concentration can increase the value of congestion awareness. In the unperturbed case, both types query the fast, inexpensive model 4 for every job under FO and SO. FO also sends every job to the more accurate but slower model 1. SO uses model 1 only when the preceding outputs disagree. Accounting for processing time therefore already reduces traffic to model 1, but SO still ignores the waiting costs generated by its chosen routes. At total arrival rate 0.8 ($\lambda_1=0.65$, $\lambda_2=0.15$), CAAR saves 34.40\% against FO and 8.16\% against SO. Savings against SO range from 5.74\% to 14.34\% across the five perturbed (random draw) cases. At total arrival rate 0.95, the base-case savings increase to 71.32\% and 14.96\%, respectively. SO's highest utilization is only 0.64 at this load. The remaining gain against SO therefore does not require a nearly saturated server. CAAR achieves these improvements by adjusting the routes as traffic increases. For example, type 2's probability of visiting model 1 is 0.506 at total arrival rate 0.2 and 0.210 at rate 0.95, while its own arrival rate remains fixed at 0.15. Thus, additional traffic from type 1 changes how CAAR routes type 2, even though type 2's demand and signal quality are unchanged.

The symmetric cases show when processing-time charges can capture much of the benefit of congestion awareness. With identical models and favorable mixing among tied policies, the main decision is how many queries to make. For one type with $b=0.5$, FO averages 2.42 queries per job. SO accounts for their processing cost and uses only one query. This is also CAAR's choice at every tested arrival rate. With two types, SO uses one query for the lower-penalty type and averages 2.42 queries for the higher-penalty type. It matches CAAR at total arrival rates 0.1, 0.3, and 0.5. Even at rate 0.9, CAAR's additional saving against SO is only 2.37\%, compared with 19.51\% against FO. Processing-time charges therefore capture much of the improvement in these cases. Waiting costs can still matter when the error penalty encourages longer routes. In the one-type case with $b=1$, SO averages 2.42 queries at every load. At arrival rate 0.9, CAAR reduces this average to 1.14 and saves 11.37\% against SO. Fixed processing-time charges can discourage excessive querying, while congestion awareness allows the query depth to adjust further as traffic increases.

\section{Conclusion}
\label{sec:conclusion}

This paper shows that designing an LLM cascade requires more than ranking models by accuracy or usage cost. Each query provides additional information for the current job, but also increases congestion at a server shared by other jobs. A model that is attractive for an individual job may therefore be a poor choice from the system's perspective. Effective cascade routing must balance the information gained from each query against both its direct cost and the additional delays it creates throughout the system.

Our CAAR framework jointly optimizes output-adaptive routing and its induced congestion. Congestion prices connect each job's information-acquisition decisions to system-wide holding costs, while column generation computes routing mixtures with a certified optimality gap. The depth bound identifies when shorter cascades suffice for a prescribed cost tolerance. The benchmark comparisons distinguish the benefits of output-dependent model selection from those of congestion-aware design. They also identify settings in which simpler policies perform nearly as well. Together, these results support designing model selection and stopping jointly with the workloads they generate.

The same perspective can guide inference tasks beyond binary classification. For multiclass prediction or regression, a specified task loss can replace the misclassification penalty, with the final decision chosen to minimize expected loss given the observed outputs. When expected losses are finite, the maintained arrival, service, and CAAR assumptions preserve the queueing formulation and convex optimization over policy mixtures. However, richer observations can make the search for these policies substantially harder. Extending the stopping and information-value results, and developing approximate pricing methods for column generation with reliable performance guarantees, are useful directions for further work.

\bibliographystyle{ormsv080}
\bibliography{CAAR_ref}

\clearpage
\appendix
\section*{Online Appendix: Proofs and Additional Results}
\addcontentsline{toc}{section}{Online Appendix: Proofs and Additional Results}
The following appendices contain the proof chains for the main-text results in
section order, followed by the one-query and two-model policies and the repeated-query extension.

\section{Proofs of Results in Section~\ref{sec:system-foundation}}
\label{apdx:proofs-1}

\begin{proof}[Proof of Lemma~\ref{lem:tree-mixture}]
For a type-$j$ job and a policy in $\Pi^{CAAR}_L$, let $R_j$ denote its entire private auxiliary randomness used by the policy. Because this randomization is independent of arrivals, service requirements, labels, and potential outputs, we may sample this auxiliary randomness in advance without changing the joint distribution of the job's routing decisions and outcomes. The horizon, feasible-history set, and action sets are finite. Conditional on $R_j=r$, the policy specifies a deterministic action after each possible observed history. These contingent decisions define a tree $T(r)\in\cT_L$: its branches correspond to model outputs, and its leaves specify terminal classifications. The realized outputs determine which path the job follows through this tree.

Define
\[
z_{jT}:=\Pp(T(R_j)=T), \qquad T\in\cT_L.
\]
Couple the original policy and the tree-mixture implementation using the same auxiliary randomness and potential outputs. Conditional on $T(R_j)=T$, both implementations follow tree $T$ pathwise, including the terminal decision. They therefore have the same joint distribution with the true label and all observed histories.

Conversely, given any distribution $(z_{jT})_{T\in\cT_L}$, use an independent uniform seed to select $T$ with these probabilities at arrival and then follow
$T$. This is an admissible CAAR policy and implements the prescribed mixture.
\end{proof}

\begin{proof}[Proof of Proposition~\ref{prop:bcmp-general}]
At arrival, pre-sample the finite vector of attributes attached to an arriving job:
\[
\mathcal M_j:=\bigl(T,Y_j,X_{j1},\ldots,X_{jM}\bigr).
\]
These attributes determine a random route $\mathsf R_j$, which is empty for immediate exit and otherwise has the form
$r=(m_1,\ldots,m_d)$ with $1\le d\le L$. Let $p_{jr}:=\Pp(\mathsf R_j=r)$, and retain only positive-length routes with $p_{jr}>0$. Independent Poisson marking gives mutually independent Poisson route streams of rate
\[
\alpha_{jr}:=\lambda_jp_{jr}
\]
for each positive-length route \cite[Section~5.2]{Kingman1993}. Sampling these job attributes at arrival leaves the distribution of the system's evolution unchanged. Each model output remains hidden from the router until the corresponding query is completed, and the sampled attributes are independent of all service requirements and the queue state immediately before the job arrives.

For each $t\in\{1,\ldots,|r|\}$, define class $(j,r,t)$ at station $m_t$.
Its transition is deterministic:
\[
(j,r,t)\longrightarrow
\begin{cases}
(j,r,t+1),&t<|r|,\\
\mathsf{exit},&t=|r|.
\end{cases}
\]
The reachable class graph is finite, open, and acyclic. Every class at physical station
$m$ receives exponential service at the same rate $\mu_m$ under FCFS.
Therefore each station is type 1 in the open BCMP theorem
\citep{BCMP1975, Kelly1979}. A job assigned route $r$ visits model $m$ exactly
$n_m(r):=\sum_{t=1}^{|r|}\mathbf1\{m_t=m\}$ times, with $n_m(r)=0$ for the empty route.
The aggregate throughput at model $m$ is therefore
\begin{align*}
\sum_{j,r}\alpha_{jr}n_m(r)
&=\sum_j\lambda_j\mathbb E[n_m(\mathsf R_j)]\\
&=\sum_j\lambda_j\sum_{T\in\cT_L}z_{jT}\mathbb E[n_m(\mathsf R_j)\mid T]\\
&=\sum_j\lambda_j\sum_{T\in\cT_L}z_{jT}v_{jTm}
=\Lambda_m(\bm z).
\end{align*}
To verify normalization and obtain the distribution of the number
of jobs at each model, let $\mathfrak C_m$ contain the route-stage
classes $\sigma=(j,r,t)$ with $\alpha_{jr}>0$ and $m_t=m$.
Each such class has throughput $\nu_\sigma:=\alpha_{jr}$.
For a station containing $n$ jobs, write
$(\sigma_1,\ldots,\sigma_n)$ for their classes in FCFS order,
starting with the job in service.
The type-1 BCMP local factor for this ordered sequence is
proportional to $\prod_{q=1}^n\nu_{\sigma_q}/\mu_m$. Hence, we have the weight for having $n$ jobs at $m$:
\begin{align*}
\sum_{(\sigma_1,\ldots,\sigma_n)\in\mathfrak C_m^n}
\prod_{q=1}^n\frac{\nu_{\sigma_q}}{\mu_m}
&=\prod_{q=1}^n \left(\sum_{\sigma\in\mathfrak C_m}\frac{\nu_\sigma}{\mu_m}\right)=\left(\frac{\sum_{\sigma\in\mathfrak C_m}\nu_\sigma}{\mu_m}\right)^n =\left(\frac{\Lambda_m(\bm z)}{\mu_m}\right)^n.
\end{align*}
Summing these weights over all possible numbers of jobs gives the geometric series $\sum_{n=0}^{\infty}(\Lambda_m(\bm z)/\mu_m)^n$.
This sum is finite, allowing the weights to be normalized into
probabilities, if and only if $\Lambda_m(\bm z)<\mu_m$. Thus, normalizing these geometric weights yields
\begin{align*}
\Pp(\mathcal N_m=n) &= \left[\sum_{n=0}^{\infty}(\Lambda_m(\bm z)/\mu_m)^n\right]^{-1}\left(\frac{\Lambda_m(\bm z)}{\mu_m}\right)^n
=\left(1-\frac{\Lambda_m}{\mu_m}\right)
\left(\frac{\Lambda_m}{\mu_m}\right)^n \mbox{ and }\\
\E[\mathcal N_m]
&= \sum_{n=0}^{\infty}n\left(1-\frac{\Lambda_m}{\mu_m}\right)\left(\frac{\Lambda_m}{\mu_m}\right)^n = \left(1-\frac{\Lambda_m}{\mu_m}\right)^{-1} \frac{\Lambda_m}{\mu_m} =\frac{\Lambda_m}{\mu_m-\Lambda_m}.
\end{align*}
This proves the stated mean-occupancy formula. Multiplying by $h$ and summing over $m$
gives the stated total holding-cost rate.

The master problem \eqref{prob:master} assigns infinite cost whenever a model's load reaches or exceeds its service rate. To justify this convention, we derive a necessary condition for stationarity.
In a stationary regime, the average arrival and departure rates
at station $m$ must coincide. Its arrival rate is $\Lambda_m(\bm z)$, while its departure rate is $\mu_m$ times the probability that the server is busy. Therefore,
\[
\Lambda_m(\bm z)
=\mu_m\Pp(\mathcal N_m>0)
\le\mu_m.
\]
Strict inequality must hold at every active station. Indeed, equality would give $\Pp(\mathcal N_m=0)=0$, whereas the reachable open network is
irreducible and its empty state must have positive probability under a
stationary distribution. Thus a mixture with
$\Lambda_m(\bm z)\ge\mu_m$ for some $m$ is not positive recurrent, which
justifies assigning it infinite steady-state cost in \eqref{prob:master}.
\end{proof}

\begin{proof}[Proof of Theorem~\ref{thm:convexity}]
The feasible set is a product of finitely dimensional probability simplices,
so it is nonempty, compact, and convex. The noncongestion term in the
objective of \eqref{prob:master} is linear in \(\bm z\). By \eqref{eq:load-tree}, each load $\Lambda_m(\bf z)$ is affine in
\(\bm z\). Recall $\Psi_m(\Lambda) = h\Lambda/(\mu_m-\Lambda)$. On \([0,\mu_m)\),
\begin{equation*}
\Psi_m''(\Lambda)
=\frac{2h\mu_m}{(\mu_m-\Lambda)^3}>0,
\end{equation*}
so \(\Psi_m\) is strictly convex in its load $\Lambda$. Its extension by \(+\infty\)
to \([0,\infty)\) is lower semicontinuous: \(\Psi_m\) is continuous on
\([0,\mu_m)\), satisfies
\(\Psi_m(\Lambda)\to+\infty\) as \(\Lambda\uparrow\mu_m\), and therefore has closed sublevel sets after the extension. Hence the objective is convex and lower semicontinuous on a compact feasible set. The zero-query tree gives a finite feasible point, so the minimum is attained.

Lastly, we show the uniqueness of the optimal load. Because a finite-cost feasible point exists, the optimal value is finite. Hence every optimal mixture has finite objective value and satisfies
\(\Lambda_m<\mu_m\) for every model. Suppose two optimal mixtures
\(\bm z^a\) and \(\bm z^b\) induce distinct load vectors. Their midpoint is feasible. Since at least one load coordinate differs, strict convexity
of the congestion cost gives
\[
C_L\left(\frac{\bm z^a+\bm z^b}{2}\right)
<\frac12C_L(\bm z^a)
+\frac12C_L(\bm z^b),
\]
contradicting optimality. Thus every optimal mixture induces the same load
vector.
\end{proof}

\begin{proof}[Proof of Theorem~\ref{thm:tree-KKT}]
Fix a stable feasible mixture $\bm z^*$. For every $m$ and every
$0\le\Lambda<\mu_m$, convexity gives the supporting inequality
\[
\Psi_m(\Lambda)-\Psi_m(\Lambda_m^*)
\ge
\Theta_m(\Lambda_m^*)(\Lambda-\Lambda_m^*).
\]
At $\Lambda_m^*>0$ this is the ordinary tangent inequality. At
$\Lambda_m^*=0$ it follows from the right derivative
$\partial_+\Psi_m(0)=h/\mu_m=\Theta_m(0)$.

Suppose first that $\bm z^*$ is optimal. Let $z_{jT}^*>0$ and
$T'\in\cT_L$. For sufficiently small $\varepsilon>0$, the mixture
$\bm z^*+\varepsilon(\bm e_{jT'}-\bm e_{jT})$ is feasible and stable. If
$\Lambda_m^*=0$, nonnegativity of the visit coefficients and
$z_{jT}^*>0$ imply $v_{jTm}=0$, so the relevant derivative is the right
directional derivative. Moreover, we can write the load shift as \(\Lambda_m(\bm z(\varepsilon)) = \Lambda_m^*+\varepsilon \lambda_j(v_{jT'm}-v_{jTm}),\) and thus $$\left.\frac{d}{d\varepsilon}
\Psi_m(\Lambda_m(\bm z(\varepsilon)))
\right|_{\varepsilon=0^+}
=
\underbrace{\Psi_m'(\Lambda_m^*)}_{\theta_m^*}
\lambda_j(v_{jT'm}-v_{jTm}).$$  
Combined with the cost part $w_{jT}$ in $A_{jT}(\bm \theta)$, optimality therefore implies
\[
\left.\frac{d}{d\varepsilon}
C_L\!\left(
\bm z^*+\varepsilon(\bm e_{jT'}-\bm e_{jT})
\right)\right|_{\varepsilon=0^+}
=
\lambda_j\bigl[
A_{jT'}(\bm\theta^*)-A_{jT}(\bm\theta^*)
\bigr]\geq 0.
\]
Because $T'$ is arbitrary and $\lambda_j>0$, every positively used tree
minimizes adjusted cost, proving necessity.

Conversely, suppose \eqref{eq:tree-KKT} holds at the induced prices
$\theta_m^*=\Theta_m(\Lambda_m^*)$. Applying the scalar supporting
inequality above to every model and using \eqref{eq:load-tree} gives, for any
stable feasible $\bm z$,
\begin{align*}
C_L(\bm z)-C_L(\bm z^*)
&\ge
\sum_{j,T}\lambda_jA_{jT}(\bm\theta^*)
(z_{jT}-z_{jT}^*)\\
&=
\sum_j\lambda_j
\left[
\sum_Tz_{jT}A_{jT}(\bm\theta^*)
-\min_{T'}A_{jT'}(\bm\theta^*)
\right]
\ge0.
\end{align*}
An unstable $\bm z$ has infinite cost, so $\bm z^*$ is globally optimal.
\end{proof}

\subsection{An optimum with few routing plans}
\label{apdx:sparse-policy}
The number of trees needed to implement an optimum is nevertheless limited. Holding the model loads fixed leaves a linear problem in the tree probabilities. This gives the following bound.

\begin{samepage}
\begin{proposition}[An optimum with few routing plans]
\label{thm:sparsity}
For every depth limit $L$, there exists an optimal CAAR mixture $\bm z^*$ satisfying
$
\bigl|\{(j,T):z_{jT}^*>0\}\bigr|\le M+J.
$
\end{proposition}
\end{samepage}

The bound applies to type--tree pairs across all
job types. For example, with ten models and three
job types, there exists an optimal policy using
at most thirteen such pairs. Since each job type
must use at least one tree, the bound allows at
most $M$ additional pairs beyond these required
choices. Thus, an optimal CAAR policy can be
implemented using relatively few routing plans,
even though each plan may involve multiple
adaptive decisions. The proof below establishes a
stronger result: any stable tree mixture can
be replaced by one satisfying this support
bound, without changing the model loads or
increasing the total cost.

\begin{proof}[Proof of Proposition~\ref{thm:sparsity}]
We establish a stronger claim for any stable feasible mixture $\widehat{\bm z}$. Let $\cC_j:=\{T:\widehat z_{jT}>0\}$ contain its trees for type $j$, and write $\widehat{\bm\Lambda}:=\bm\Lambda(\widehat{\bm z})$. Hold these loads fixed and solve
\[
\begin{aligned}
\min_{\bm z\ge0}\quad&
\sum_j\lambda_j\sum_{T\in\cC_j}z_{jT}\omega_{jT}\\
\text{s.t.}\quad&
\sum_{T\in\cC_j}z_{jT}=1,\qquad j\in[J],\\
&\sum_j\lambda_j\sum_{T\in\cC_j}z_{jT}v_{jTm}
=\widehat\Lambda_m,\qquad m\in[M].
\end{aligned}
\]
The original mixture is feasible, and the feasible set is a bounded polytope. The linear objective has an optimal extreme point $\widetilde{\bm z}$. There are at most $M+J$ independent equality constraints, so this extreme point has at most $M+J$ positive weights \citep{BertsimasTsitsiklis1997}.

The new mixture uses only trees in the original mixture and preserves all aggregate model loads. It therefore preserves stability and holding costs. Aggregate query fees also remain unchanged because they equal $\sum_m c_m\widehat\Lambda_m$. The error-penalty rate cannot increase, so $C_L(\widetilde{\bm z})\le C_L(\widehat{\bm z})$. Applying this construction to an optimal mixture proves the proposition.
\end{proof}

\section{Proofs of Results in Section~\ref{sec:policy-structure}}
\label{apdx:proofs-2}

\begin{proof}[Proof of Lemma~\ref{lem:posterior-sufficiency}]
Represent all private randomization by a seed $U_j$ uniformly distributed on $[0,1]$, independent of the label and potential outputs.
Let $M_{j,s}$ denote the model queried at stage $s$, and define
\[
H_{j,t}^{\mathrm{obs}}
:=\bigl((M_{j,s},X_{j,M_{j,s}})\bigr)_{s=1}^t,
\qquad H_{j,0}^{\mathrm{obs}}:=\varnothing.
\]
If the job stops before stage $t$, fill the remaining entries with $\dagger$.
The information generated by the seed and observed history is
$\mathcal H_{j,t}:=\sigma(U_j,H_{j,t}^{\mathrm{obs}})$, with
$\mathcal H_{j,0}=\sigma(U_j)$.
Write the action at stage $s$ as the measurable function
$\varphi_s(u,\mathfrak h_{s-1})$. For a reachable
$\mathfrak h_t=((m_1,x_1),\ldots,(m_t,x_t))$, let
\[
\mathsf H_t
:=\{H_{j,t}^{\mathrm{obs}}=\mathfrak h_t\}.
\]
For a Borel
$B\subseteq[0,1]$, define the label-independent selection factor
\[
a_B(\mathfrak h_t)
:=\int_B\prod_{s=1}^t
\mathbf1\{\varphi_s(u,\mathfrak h_{s-1})=m_s\}\,du.
\]
Independence of the seed and conditional signal independence imply, for
$y\in\{A,B\}$,
\begin{align}
&\Pp\bigl(
Y_j=y,\ U_j\in B,\ \mathsf H_t
\bigr) =
\Pp(Y_j=y)a_B(\mathfrak h_t)
\prod_{s=1}^t\ell_{jm_s}^y(x_s).
\label{eq:posterior-selection-factor}
\end{align}
For every \(B\) for which the conditioning event has positive probability,
the common factor \(a_B(\mathfrak h_t)\) cancels in Bayes' rule. Because
\eqref{eq:posterior-selection-factor} holds for every Borel \(B\), these ratios
define the same version of the posterior conditional on
\(\mathcal H_{j,t}\). This gives the stated posterior formula almost surely
on $\mathsf H_t$.

We next show the mutual independence of unused models. For $S\subseteq[M]\setminus\{m_1,\ldots,m_t\}$, the corresponding joint
factorization is
\begin{align*}
&\Pp\bigl(
Y_j=y,\ U_j\in B,\ \mathsf H_t,\
(X_{jm})_{m\in S}=x_S
\bigr) =
\Pp(Y_j=y)a_B(\mathfrak h_t)
\left(\prod_{s=1}^t\ell_{jm_s}^y(x_s)\right)
\left(\prod_{m\in S}\ell_{jm}^y(x_m)\right).
\end{align*}
On conditioning events of positive probability, division by
\eqref{eq:posterior-selection-factor} yields the conditional independence
and likelihoods of the unused signals stated in the lemma. On
$\mathsf H_t$, the conditional law of total probability then gives the stated predictive probability
\begin{align*}
\Pp(X_{jm}=x\mid\mathcal H_{j,t})
&=\sum_{y\in\{A,B\}}
\Pp(X_{jm}=x\mid Y_j=y,\mathcal H_{j,t})
\Pp(Y_j=y\mid\mathcal H_{j,t})\\
&=p_j(\mathfrak h_t)\ell_{jm}^A(x)
+\bigl(1-p_j(\mathfrak h_t)\bigr)\ell_{jm}^B(x).
\end{align*}
One more application of Bayes' rule gives $\tau_{jm,x}$ on the same event.
\end{proof}

\begin{proof}[Proof of Lemma~\ref{lem:MAP-WLOG}]
Fix type $j$ and let $\mathcal I$ be the information generated by the seed and
the observed signal history when routing stops, before the terminal label is
selected. Write $p_{\mathcal I}:=\Pp(Y_j=A\mid\mathcal I)$. If a terminal rule
selects $A$ with conditional probability $\pi_A\in[0,1]$, then
\begin{align*}
\Pp(\widehat Y\ne Y_j\mid\mathcal I)
&=\pi_A(1-p_{\mathcal I})+(1-\pi_A)p_{\mathcal I} \ge \min\{p_{\mathcal I},1-p_{\mathcal I}\},
\end{align*}
with equality for a MAP label. Replacing only the terminal rule leaves the
realized route, query fees, service requirements, and model loads unchanged.
Taking expectations proves the result.
\end{proof}

\begin{proof}[Proof of Proposition~\ref{prop:Bellman}]
By Lemma \ref{lem:posterior-sufficiency}, at every reachable history the actual continuation problem is exactly \eqref{eq:W-definition}, with \(p\) equal to the realized posterior. Thus
\((p,S,r)\) is a sufficient continuation state. If \(r=0\), MAP stopping
gives $W_{j,S,0}(p;\bm\chi)=b_j\min\{p,1-p\}=g_j(p).$

Now let \(r\ge1\). The first action partitions the feasible policies into
immediate stopping and one class for each \(m\in S\). Conditional on querying
\(m\), output \(x\) has probability \(s_{jm,x}(p)\), and the continuation
state is thus $\bigl(\tau_{jm,x}(p),S\setminus\{m\},r-1\bigr).$
Because the finite continuation-policy sets attain their minima and the two output branches can be selected independently,
\begin{align*}
&\min_{\pi\in\Pi^{c}_{j,S,r}:\,\pi\text{ queries }m\text{ first}}
 \mathbb E_{j,p,S}^{\pi}
 \left[
 \sum_{t=1}^{N^\pi}\chi_{M_t^\pi}
 +b_j\mathbf 1\{\widehat Y^\pi\ne Y_j\}
 \right]\\
&\qquad=
 \chi_m+\sum_{x\in\{A,B\}}s_{jm,x}(p)
 W_{j,S\setminus\{m\},r-1}
 \bigl(\tau_{jm,x}(p);\bm\chi\bigr)
 =G_{jm,S,r}(p;\bm\chi).
\end{align*}
Taking the minimum over these first-action classes proves
\eqref{eq:Bellman} and the optimal action rule.

Finally, by Lemma \ref{lem:tree-mixture}, randomization cannot improve upon the
finite set \(\cT_L\) of deterministic trees. For each \(T\in\cT_L\),
\[
\mathbb E\left[
\sum_{t=1}^{N^T}\chi_{M_t^T}
+b_j\mathbf 1\{\widehat Y^T\ne Y_j\}
\right]
=\sum_{m=1}^M\chi_m v_{jTm}+b_je_{jT}.
\]
Minimizing this identity over \(T\in\cT_L\) proves
\eqref{eq:root-tree-equivalence}.
\end{proof}

\begin{proof}[Proof of Corollary~\ref{cor:equilibrium-policy}]
By Proposition \ref{prop:Bellman}, a tree has minimum adjusted cost if and only if its
action is Bellman optimal at every node that the tree can reach with positive probability. Specifications after a stopping node do not affect its cost or visit vector. The claim now follows in both directions from
Theorem \ref{thm:tree-KKT}.
\end{proof}

\begin{proof}[Proof of Lemma~\ref{lem:posterior-geometry}]
Fix $j,S,r,\bm\chi$ and a complete deterministic continuation policy $\pi$,
including its terminal labels. Conditional on $Y_j=y$, the future outputs
follow the fixed type-$j$ likelihoods and the actions prescribed by $\pi$.
Its conditional expected query and error cost is therefore a finite
nonnegative quantity $W^{\pi,y}(\bm\chi)$, independent of $p$. That is, the law of total expectation gives
\begin{align*}
W^\pi(p;\bm\chi) &= 
p\mathbb E\!\left[
\sum_{t=1}^{N^\pi}\chi_{M_t^\pi}
+b_j\mathbf1\{\widehat Y^\pi\ne A\}
\,\middle|\,Y_j=A 
\right] + (1-p) \mathbb E\!\left[
\sum_{t=1}^{N^\pi}\chi_{M_t^\pi}
+b_j\mathbf1\{\widehat Y^\pi\ne B\}
\,\middle|\,Y_j=B 
\right]\\
&=pW^{\pi,A}(\bm\chi)+(1-p)W^{\pi,B}(\bm\chi).
\end{align*}
There are finitely many complete deterministic continuation policies,
and randomization cannot lower their minimum cost. Hence
$W_{j,S,r}(\cdot;\bm\chi)$ is the pointwise minimum of finitely many
affine functions. It is consequently continuous, concave, and piecewise
affine with finitely many pieces in $p$.

All query and error costs are nonnegative. Immediate MAP stopping is
feasible and has cost $g_j(p)$, proving the value bounds. At $p=0$ or
$p=1$, the label is known and immediate stopping has zero cost, so
nonnegativity gives the endpoint equalities.

For $r\ge1$ and $m\in S$, $G_{jm,S,r}(p;\bm\chi)$ is the minimum cost over
the finite subset of complete policies that query $m$ first. The same
affine-cost argument gives its continuity, concavity, and finite
piecewise-affine representation in $p$.
\end{proof}

\begin{proof}[Proof of Proposition~\ref{prop:stopping-threshold}]
On the interval \([0,1/2]\), the stopping cost is the affine function
\(g_j(p)=b_jp\). Define
\begin{equation*}
d_-(p):=b_jp-W_{j,S,r}(p;\bm\chi),
\qquad p\in[0,1/2].
\end{equation*}
By Lemma \ref{lem:posterior-geometry}, \(d_-\) is nonnegative and convex,
with \(d_-(0)=0\).

The zero set of a nonnegative convex function is convex. Indeed, if
\(d_-(p_1)=d_-(p_2)=0\), then for every \(t\in[0,1]\), convexity gives
\[
0\le d_-(tp_1+(1-t)p_2)
\le td_-(p_1)+(1-t)d_-(p_2)=0.
\]
Thus the zero set of \(d_-\) is an interval containing \(0\). Since it is
closed by continuity, it has the form \([0,p^-]\) for some
\(p^-\in[0,1/2]\).

On \([1/2,1]\), the stopping cost is the affine function
\(g_j(p)=b_j(1-p)\). Define
\begin{equation*}
d_+(p):=b_j(1-p)-W_{j,S,r}(p;\bm\chi).
\end{equation*}
The same argument shows that \(d_+\) is nonnegative and convex, with
\(d_+(1)=0\). Its zero set is therefore a closed interval of the form
\([p^+,1]\), where \(p^+\in[1/2,1]\). Combining the two halves gives
the stated stopping set. On $(p^-,p^+)$, $W_{j,S,r}(p;\bm\chi)<g_j(p)$, so querying is strictly preferable to stopping.
If $p^-<p^+$, then $r\ge1$ and the Bellman equation gives
$\min_{m\in S}G_{jm,S,r}(p;\bm\chi)=W_{j,S,r}(p;\bm\chi)$ on this interval.
Both sides are continuous, so taking limits from within the interval gives equality with $g_j$ at either cutoff.
If $p^-=p^+=1/2$, stopping is optimal everywhere, and a tie with querying is not required.

It remains to prove the explicit tail bounds. Any action that queries a model
incurs an immediate adjusted cost of at least \(\underline\chi_S\), followed
by a nonnegative future cost. Hence every continuation action costs at least
\(\underline\chi_S\). If
\(p\le1/2\) and \(b_jp\le\underline\chi_S\), stopping is weakly cheaper
than every continuation action. Therefore
\[
\left[0,\min\left\{\frac12,
\frac{\underline\chi_S}{b_j}\right\}\right]
\subseteq\mathcal Z_{j,S,r},
\]
which gives lower bound. The upper bound follows identically, using
\(b_j(1-p)\le\underline\chi_S\).
\end{proof}

\begin{proof}[Proof of Corollary~\ref{cor:cost-stopping}]
For any fixed complete continuation policy, increasing adjusted costs
componentwise on $S$ weakly increases its expected query cost and leaves
its expected error cost unchanged. The feasible policy set is the same,
so taking minima gives $W_{j,S,r}(p;\bm\chi')\ge W_{j,S,r}(p;\bm\chi).$

If $p\in\mathcal Z_{j,S,r}(\bm\chi)$, the right-hand side is $g_j(p)$.
Immediate stopping is still feasible under $\bm\chi'$, so the left-hand
side cannot exceed $g_j(p)$. Equality then follows, proving the stopping-set inclusion.
\end{proof}

\begin{proof}[Proof of Lemma~\ref{lem:W-lipschitz}]
Fix \(j,S,r,\bm\chi\) and write \(W(p):=W_{j,S,r}(p;\bm\chi)\).
By Lemma \ref{lem:posterior-geometry}, \(W\) is concave, satisfies
\(W(0)=W(1)=0\), and obeys $0\leq W(p) \leq b_j\min\{p,1-p\}$. Fix \(0\le p<q\le1\).  Suppose first that \(p>0\).  Since
\(p=(p/q)q+(1-p/q)0\), concavity and \(W(0)=0\) give $W(p)\ge \frac{p}{q}W(q).$
Therefore, we have
\[
 W(q)-W(p)
 \le \frac{q-p}{p}W(p)
 \le b_j(q-p),
\]
where the last inequality follows from \(W(p)\le b_jp\).  The same conclusion
holds when \(p=0\), because then
\(W(q)-W(0)=W(q)\le b_jq\).

For the reverse difference, suppose first that \(q<1\).  Since $q= [(1-q)p
+(q-p)1]/(1-p),$ concavity and \(W(1)=0\) give
\[
 W(q)\ge \frac{1-q}{1-p}W(p).
\]
Hence,
\[
 W(p)-W(q)
 \le \frac{q-p}{1-q}W(q)
 \le b_j(q-p),
\]
using \(W(q)\le b_j(1-q)\).  When \(q=1\), the same inequality follows from
\(W(p)\le b_j(1-p)\). Combining the two directions proves the bound.
\end{proof}

\begin{proof}[Proof of Theorem~\ref{thm:statewise-information-frontier}]
Fix \(m\in S\) and suppose that querying \(m\) is Bellman optimal.  Define
\[
 \widetilde W(q)
 :=
 W_{j,S\setminus\{m\},r-1}(q;\bm\chi).
\]
The policy that ignores model \(m\) and follows an optimal continuation using
at most \(r-1\) models from \(S\setminus\{m\}\) is feasible at the original
state.  Therefore Bellman optimality of query \(m\) implies
\begin{align*}
 \chi_m
 &+
 \sum_{x\in\{A,B\}}
 s_{jm,x}(p)
 \widetilde W\bigl(\tau_{jm,x}(p)\bigr) =W_{j,S,r}(p;\bm\chi)
 \le \widetilde W(p).
\end{align*}
Let \(P'\) denote the posterior after observing model \(m\)'s output.  By
Lemma \ref{lem:W-lipschitz} and the optimality $\chi_m+\mathbb E[\widetilde W(P')] \leq \widetilde W(p)$, we have
\begin{align}
 \chi_m
 &\le
 \widetilde W(p)-\mathbb E[\widetilde W(P')]
 \le
 \mathbb E\left[
 \left|\widetilde W(p)-\widetilde W(P')\right|
 \right]
 \le
 b_j\mathbb E|P'-p|.
\label{eq:frontier-lipschitz-step}
\end{align}
Positive orientation gives
\(\tau_{jm,A}(p)\ge p\ge\tau_{jm,B}(p)\).  Directly from
\eqref{eq:predictive}, \eqref{eq:update} and $d_{jm}=q_{jm}^A+q_{jm}^B-1$,
\begin{align*}
 s_{jm,A}(p)
 \left(\tau_{jm,A}(p)-p\right)
 &=p(1-p)d_{jm},
\\
 s_{jm,B}(p)
 \left(p-\tau_{jm,B}(p)\right)
 &=p(1-p)d_{jm}.
\end{align*}
Adding the two identities yields $\mathbb E|P'-p| =2p(1-p)d_{jm}.$
Combining with \eqref{eq:frontier-lipschitz-step} thus gives
\begin{equation}
 \chi_m\le2b_jp(1-p)d_{jm}.
\label{eq:frontier-info-cap}
\end{equation}

A second bound follows from immediate stopping.  The future value after a query
is nonnegative, while Bellman optimality and feasibility of stopping imply
\[
 \chi_m
 \le G_{jm,S,r}(p;\bm\chi)
 =W_{j,S,r}(p;\bm\chi)
 \le g_j(p)
 =b_j\min\{p,1-p\}.
\]
Combining this inequality with \eqref{eq:frontier-info-cap} proves \eqref{eq:statewise-frontier}. %
\end{proof}

\begin{proof}[Proof of Corollary~\ref{cor:exogenous-frontier}]
Suppose type \(j\) queries model \(m\) with positive probability under an optimal mixture.  Then \(\Lambda_m^*>0\), and some positively weighted tree queries \(m\) at a history that the tree reaches with positive probability.
By Corollary \ref{cor:equilibrium-policy}, querying \(m\) is Bellman optimal at that history under the equilibrium adjusted cost
$\chi_m^*=c_m+\Theta_m(\Lambda_m^*).$ Applying Theorem \ref{thm:statewise-information-frontier} and using
\(2p(1-p)\le1/2\) gives
\begin{equation*}
 \chi_m^*
 \le \frac{b_jd_{jm}}{2}.
\end{equation*}
Because \(\Lambda_m^*>0\) and \(\Theta_m\) is strictly increasing, we have
$\chi_m^* >c_m+\Theta_m(0) =\chi_m^{\mathrm{idle}}.$
Hence positive type-\(j\) use requires
\(\chi_m^{\mathrm{idle}}<b_jd_{jm}/2\), proving the zero-visit claim in
part~(i).  Any specification of model \(m\) in a tree but not in the set can be replaced by stopping without changing the tree's cost or visit vector,
so restricting the tree space is without loss.  Part~(ii) follows by taking
the union across types. Next, if \(\Lambda_m^*>0\), the same argument for a type that uses model \(m\)
gives
\[
 c_m+
 \frac{h\mu_m}{(\mu_m-\Lambda_m^*)^2}
 =\chi_m^*
 \le\overline v_m.
\]
Rearranging the preceding inequality yields
\eqref{eq:information-load-cap}. In particular,
\(\overline v_m>c_m+h/\mu_m=\chi_m^{\mathrm{idle}}\). Its contrapositive gives the inactivity
statement when \(\overline v_m\le\chi_m^{\mathrm{idle}}\). This proves Part (iii).
\end{proof}

\section{Proofs of Results in Section~\ref{sec:column-generation}}
\label{apdx:proofs-3}

\begin{proof}[Proof of Proposition~\ref{thm:cg-certificate}]
Fix $K$ and retain the immediate-stopping tree in each generated set. The restricted-master feasible set is a finite product of simplices. Its objective is lower semicontinuous and convex, and immediate stopping has finite cost. Hence a minimum is attained, and every minimizer is stable.

Let $\bm z$ be a restricted-master minimizer with cost $U$ and prices $\bm\theta$. A small shift of type-$j$ probability from an active tree $T$ to any $T'\in\cC_j$ remains stable. Optimality gives $
\lambda_j\bigl[A_{jT'}(\bm\theta)-A_{jT}(\bm\theta)\bigr]\ge0.$
Thus, every active tree has adjusted cost $\nu_j$, and
$\sum_{T\in\cC_j}z_{jT}A_{jT}(\bm\theta)=\nu_j$.
Also, $\phi_j(\bm\theta)\le\nu_j$ because $\cC_j\subseteq\cT_{j,K}$.

For any stable mixture $\widetilde{\bm z}$ over the full screend set $\cT_{j,K}$, convexity of each $\Psi_m$ implies
\[
\Psi_m(\Lambda_m(\widetilde{\bm z}))
\ge
\Psi_m(\Lambda_m(\bm z))
+\theta_m\bigl[\Lambda_m(\widetilde{\bm z})-\Lambda_m(\bm z)\bigr].
\]
At zero load, the same inequality uses the right derivative $\Theta_m(0)$.
With \(\frac{\partial C_K}{\partial z_{jT}}(\bm z)
=\lambda_j A_{jT}(\bm\theta)\), adding the linear query and error costs yields
\[
\begin{aligned}
C_K(\widetilde{\bm z})
&\ge U+\sum_j\lambda_j
\left[\sum_{T\in\cT_{j,K}}\widetilde z_{jT}A_{jT}(\bm\theta)-\nu_j\right] \ge U+\sum_j\lambda_j\bigl[\phi_j(\bm\theta)-\nu_j\bigr]
=LB.
\end{aligned}
\]
For unstable mixtures, this bound is immediate because the cost is infinite.
Corollary~\ref{cor:exogenous-frontier} ensures that minimizing over the screened sets gives $C_K^*$. Hence $LB\le C_K^*\le U$.

If every pricing equality holds, then $LB=U$, proving global optimality. Conversely, suppose $\phi_j(\bm\theta)<\nu_j$ for some type. Let $T_j$ attain $\phi_j(\bm\theta)$ and shift a positive amount $\eta$ of probability from an active tree to $T_j$. For sufficiently small $\eta$, stability is preserved and we have $
\lambda_j\bigl[\phi_j(\bm\theta)-\nu_j\bigr]<0.$
This gives a strictly lower cost, contradicting global optimality.
Such a tree cannot belong to $\cC_j$. Every nonterminal iteration therefore adds at least one previously absent tree. Since the sets $\cT_{j,K}$ are finite and generated columns are retained, only finitely many additions are possible. At termination all pricing equalities hold, so the value is $C_K^*$.
\end{proof}

The same argument gives a bound when the restricted master is solved only approximately. For any stable feasible mixture $\widehat{\bm z}$, use its induced prices and set
$\overline A_j:=\sum_T\widehat z_{jT}A_{jT}(\bm\theta)$. Exact tree pricing then gives
\[
C_K(\widehat{\bm z})-\sum_j\lambda_j
\bigl[\overline A_j-\phi_j(\bm\theta)\bigr]
\le C_K^*\le C_K(\widehat{\bm z}).
\]
Thus the certificate can be checked without assuming that a numerical master solution satisfies the exact optimality conditions. Each $\phi_j$ here must be the true pricing minimum, or be replaced by a valid lower bound on that minimum.

\section{Proofs of Results in Section~\ref{sec:effective-depth}}
\label{apdx:proofs-5}

\paragraph{Signal overlap and the adaptive frontier.}
Interior signal probabilities give $\beta_{jm}>0$, and Cauchy--Schwarz
gives $\beta_{jm}\le1$. Equality requires the two label-conditional
distributions to coincide, or $q_{jm}^A=1-q_{jm}^B$, which
Assumption \ref{ass:orientation} excludes. Thus $0<\beta_{jm}<1$ and
$\mathcal I_{jm}>0$.

For a deterministic tree $T$, let $\mathfrak H_{j,r}(T)$ contain the
histories at which a type-$j$ job completes $r$ queries without having
stopped earlier. Write $P_{j,y}^T(\mathfrak h)$ for the probability of
reaching $\mathfrak h$ conditional on $Y_j=y$. At $r=0$ the frontier
contains the empty history.

\begin{lemma}[Overlap and Bayes error]
\label{lem:adaptive-frontier-affinity}
For every type $j$, deterministic tree $T$, depth $r\in\{0,\ldots,M\}$,
and subset $\mathfrak F\subseteq\mathfrak H_{j,r}(T)$,
\begin{align*}
\sum_{\mathfrak h\in\mathfrak F}
\sqrt{P_{j,A}^T(\mathfrak h)P_{j,B}^T(\mathfrak h)}
&\le e^{-\mathcal I_j(r)},\\
\sum_{\mathfrak h\in\mathfrak F}
\min\{\rho_jP_{j,A}^T(\mathfrak h),
(1-\rho_j)P_{j,B}^T(\mathfrak h)\}
&\le\sqrt{\rho_j(1-\rho_j)}e^{-\mathcal I_j(r)}.
\end{align*}
\end{lemma}

\begin{proof}[Proof of Lemma~\ref{lem:adaptive-frontier-affinity}]
Fix an available set $S$ and $0\le r\le |S|$.
Consider any deterministic continuation tree that queries distinct
models from $S$ along each path. For its histories reached after
$r$ additional queries, excluding paths that stop earlier, measure
overlap by summing the geometric mean of each history's probabilities
under labels $A$ and $B$. We show by induction on $r$, for every such
set and tree, that this overlap is at most the product of the $r$
largest coefficients $\{\beta_{jm}:m\in S\}$.

For $r=0$, the only history is the empty history. Its probability
is one under either label, so the overlap and the empty product
both equal one.

For $r\ge1$, suppose the claim holds for $r-1$ additional queries
for every available set and continuation tree. If the tree stops
at the current node, no history reaches $r$ additional queries,
so the overlap is zero. Otherwise, the tree queries a model
$m\in S$. After its output $x$, the remaining tree uses models
in $S\setminus\{m\}$. Conditional independence ensures that
the unqueried outputs keep their original label-conditional
distributions. The choice of $m$ is fixed by the current history
and introduces no additional likelihood factor.

By the induction hypothesis, the overlap of the histories reached
after $r-1$ further queries in either remaining tree is at most
the product of the $r-1$ largest coefficients in $S\setminus\{m\}$.
For the branch corresponding to output $x$, each remaining-history
contribution is multiplied by
$\sqrt{\ell_{jm}^A(x)\ell_{jm}^B(x)}$.
Applying the same bound to both branches and summing therefore
bounds the overlap after $r$ queries by
\[
\sum_{x\in\{A,B\}}
\sqrt{\ell_{jm}^A(x)\ell_{jm}^B(x)}
\prod_{k=1}^{r-1}\beta_{j,S\setminus\{m\},[k]}
=
\beta_{jm}\prod_{k=1}^{r-1}\beta_{j,S\setminus\{m\},[k]} \leq \prod_{k=1}^{r}\beta_{j,[M],[k]}
=e^{-\mathcal I_j(r)},
\]
where bracketed indices indicate decreasing order.
This is a product of $r$ coefficients from $S$, so it cannot
exceed the product of the $r$ largest. This completes the induction.

Taking $S=[M]$ proves the first inequality for all histories
in $\mathfrak H_{j,r}(T)$, and nonnegativity extends it to any subset.
Applying $\min\{a,b\}\le\sqrt{ab}$ at each history gives the second.
\end{proof}

\begin{proof}[Proof of Theorem~\ref{thm:adaptive-depth-truncation}]
Take an optimal depth-$L$ tree mixture. Replace each tree $T$ by
$T^{[K]}$, which follows $T$ until it stops or completes $K$ queries.
At a depth-$K$ history where $T$ would continue, $T^{[K]}$ stops and
uses MAP. Keeping the original mixture probabilities gives a feasible
depth-$K$ policy.

On every realization, truncation only removes future queries.
Thus each visit probability and model load weakly decreases.
The original optimum is stable, so the truncated policy remains stable.
Query fees weakly decrease, as does holding cost because
$\Psi_m(\Lambda)$ is increasing for $0\le\Lambda<\mu_m$.

For type $j$, the histories where truncation occurs form a subset
$\mathfrak F\subseteq\mathfrak H_{j,K}(T)$. At such a history
$\mathfrak h$, the joint error mass of immediate MAP stopping is
$\min\{\rho_jP_{j,A}^T(\mathfrak h),
(1-\rho_j)P_{j,B}^T(\mathfrak h)\}$.
Let \(\eta_j^T(\mathfrak h)\ge0\) denote the joint probability that the job reaches \(\mathfrak h\) and is eventually misclassified under the original continuation, so
Lemma \ref{lem:adaptive-frontier-affinity} gives
\[
\begin{aligned}
e_j(T^{[K]})-e_j(T)
&=\sum_{\mathfrak h\in\mathfrak F}
\Bigl[
\min\{\rho_jP_{j,A}^T(\mathfrak h),
(1-\rho_j)P_{j,B}^T(\mathfrak h)\}
-\eta_j^T(\mathfrak h)
\Bigr]\\
&\le
\sum_{\mathfrak h\in\mathfrak F}
\min\{\rho_jP_{j,A}^T(\mathfrak h),
(1-\rho_j)P_{j,B}^T(\mathfrak h)\}\le
\sqrt{\rho_j(1-\rho_j)}e^{-\mathcal I_j(K)}.
\end{aligned}
\]
Averaging over trees preserves this bound. Multiplying by $\lambda_jb_j$
and summing over types bounds the increase in total cost by
$\mathcal E(K)$. Optimality of $C_K^*$ then gives the upper bound in
the theorem. Inclusion of the depth-$K$ class in the depth-$L$ class gives the lower bound. Setting the larger depth to $L+1$ gives the marginal-depth inequality
in the main text. A simpler type-wise bound follows from
$\mathcal I_j(K)\ge K\underline{\mathcal I}_j$, where
$\underline{\mathcal I}_j:=\min_m\mathcal I_{jm}$:
$C_K^*-C_L^*\le\sum_j\lambda_jB_j e^{-K\underline{\mathcal I}_j}$.
\end{proof}

\begin{proof}[Proof of Corollary~\ref{cor:adaptive-epsilon-depth}]
The defining set for $K_\varepsilon$ is nonempty because
$\mathcal E(M)=0$. If $K_\varepsilon=M$, the cost gap is zero.
Otherwise $K_\varepsilon<M$ and
$\mathcal E(K_\varepsilon)\le\varepsilon$, so
Theorem \ref{thm:adaptive-depth-truncation} with $L=M$ proves
the claimed guarantee. This includes $K_\varepsilon=0$.

For the explicit depth, $\mathcal I_j(K)\ge K\underline{\mathcal I}$
implies $\mathcal E(K)\le B_\Sigma e^{-K\underline{\mathcal I}}$ for
$K<M$. If $\widehat K_\varepsilon=M$, the gap is again zero.
If $\widehat K_\varepsilon=0<M$, its definition implies
$\varepsilon\ge B_\Sigma=\mathcal E(0)$.
In the remaining case $0<\widehat K_\varepsilon<M$, the ceiling ensures
$B_\Sigma e^{-\widehat K_\varepsilon\underline{\mathcal I}}\le\varepsilon$.
Applying the truncation theorem proves the claim in both cases.
\end{proof}

\section{Proof of the Fixed-Cascade Comparison}
\label{apdx:fixed-cascades}

\begin{proof}[Proof of Proposition~\ref{prop:fixed-ranking-gap}]
Fix $h>0$. Consider one job type with prior $\rho=1/2$, arrival rate $\lambda=1$, and at most two queries. All three models have fee $c_m=1$ and service rate $\mu_m=2$. For $0<\varepsilon\le1/8$, let $\eta=\varepsilon^2$, $b=\varepsilon^{-2}$, and set
\[
(q_1^A,q_1^B)=(1-\varepsilon,1-\varepsilon),\qquad
(q_2^A,q_2^B)=(1-\eta,1-\varepsilon),\qquad
(q_3^A,q_3^B)=(1-\varepsilon,1-\eta).
\]
The outputs are conditionally independent given the label. Model 1 is the screener. Models 2 and 3 are specialists with complementary accuracies. Every model's load is at most $\lambda=1<2$, so every admissible policy is stable.

Under the balanced prior, the Bayes error of an observation $Z$ is
\[
\frac12\sum_z\min\{\Pp(Z=z\mid Y=A),\Pp(Z=z\mid Y=B)\}.
\]
Either specialist alone has error $e_{\mathrm{sp}}=(\varepsilon+\eta)/2$. For the pair of models 1 and 2, the MAP decision always follows model 2. On disagreements, the relevant likelihood comparisons are
$(1-\varepsilon)\eta<\varepsilon(1-\varepsilon)$ and
$\varepsilon(1-\eta)>\varepsilon(1-\varepsilon)$.
Agreement outcomes favor their common label. The pair error therefore remains $e_{\mathrm{sp}}$. Exchanging $A$ and $B$ gives the same conclusion for models 1 and 3. For models 2 and 3, summing the smaller likelihoods over the four outputs gives $2\varepsilon\eta+\varepsilon(1-\eta)+\eta(1-\varepsilon)
=\varepsilon+\eta,$ implying the Bayes error $e_{\mathrm{sp}}$. Thus, every model pair has Bayes error $e_{\mathrm{sp}}$.

Now take any FR or MR policy $\pi$, including arbitrary randomized stopping. Let $\alpha$ be its probability of making an initial query. Conditional on its chosen ranking and private randomization, an observer of the first two outputs in that ranking can reproduce the policy's terminal decision, including any early stopping. The private randomization and the initial-query decision are independent of the label and outputs. Conditional on querying, the policy's error is therefore at least $e_{\mathrm{sp}}$, and it pays at least one query fee. Conditional on stopping immediately, its error is at least $1/2$. Holding costs are nonnegative. Let \(C_2(\pi)\) denote the total steady-state cost per unit time under policy \(\pi\), when each job may make at most two queries. Its total cost satisfies
\[
C_2(\pi)\ge (1-\alpha)\frac b2+\alpha(1+b e_{\mathrm{sp}})
\ge\min\left\{\frac b2,\,1+b e_{\mathrm{sp}}\right\}.
\]
For the chosen parameters,
\[
1+b e_{\mathrm{sp}}=\frac32+\frac1{2\varepsilon},\qquad
\frac b2-(1+b e_{\mathrm{sp}})
=\frac{1-\varepsilon-3\varepsilon^2}{2\varepsilon^2}>0.
\]
This establishes $C_{\mathrm{MR}}^*\ge3/2+1/(2\varepsilon)$, and the same bound applies to FR. It includes the option of immediate stopping.

For an upper bound on the CAAR optimum, query model 1, then model 2 after output $A$ or model 3 after output $B$, and finally use MAP. Each first output has probability $1/2$. After $A$, the posterior is $1-\varepsilon$, and querying model 2 gives conditional error from MAP:
\[
\min\{(1-\varepsilon)(1-\eta),\varepsilon^2\}
+\min\{(1-\varepsilon)\eta,\varepsilon(1-\varepsilon)\}
=\varepsilon^2+(1-\varepsilon)\eta.
\]
The same error follows after observing $B$ by symmetry. Hence this feasible tree has error $e_{\mathrm{br}}=\varepsilon^2+(1-\varepsilon)\eta$, pays two query fees, and induces loads $(1,1/2,1/2)$. Its holding cost is
\[
h\left(\frac{1}{2-1}
+\frac{1/2}{2-1/2}
+\frac{1/2}{2-1/2}\right)=\frac{5h}{3}.
\]
Its total cost is therefore $2+b e_{\mathrm{br}}+5h/3=4-\varepsilon+5h/3$. Since this bounds $C_2^*$ from above and FR is a subclass of MR,
\[
\frac{C_{\mathrm{FR}}^*}{C_2^*}
\ge\frac{C_{\mathrm{MR}}^*}{C_2^*}
\ge\frac{3/2+1/(2\varepsilon)}{4-\varepsilon+5h/3}
\longrightarrow\infty
\qquad\text{as }\varepsilon\downarrow0.
\]
The denominator is the cost of a feasible CAAR policy. The argument does not require this policy to be optimal for the full-cost objective.
\end{proof}

The dependence on the error penalty is necessary for an unbounded ratio when immediate stopping is available. To see this, consider any one-type instance with $p_0=\min\{\rho,1-\rho\}>0$ and $c_{\min}=\min_m c_m>0$. A CAAR policy that makes an initial query with probability $\alpha$ has total cost at least
\[
\lambda\big[(1-\alpha)b p_0+\alpha c_{\min}\big]
\ge\lambda\min\{b p_0,c_{\min}\}.
\]
Immediate stopping belongs to both ranking classes and costs $\lambda b p_0$. Consequently,
\[
1\le\frac{C_{\mathrm{MR}}^*}{C_K^*}
\le\frac{C_{\mathrm{FR}}^*}{C_K^*}
\le\max\left\{1,\frac{b p_0}{c_{\min}}\right\}.
\]
This bound concerns total cost, including holding costs. It shows why the construction increases $b/c_{\min}$ while keeping arrival rates, fees, and service rates fixed.

\section{One-Query and Two-Model Policies}
\label{apdx:policy-special-cases}

\subsection{One query: model selection and water-filling}
\label{subsec:L1}

With at most one query per job ($L=1$), model selection reduces to comparing each query with immediate classification. The system may still use many models across jobs. Symmetric accuracies give a simple index for this comparison. If accuracies and error penalties are also common across types, one scalar determines the optimal loads. 

We first evaluate a query at fixed adjusted costs. Because no further query is allowed, the job uses the output for MAP classification and then leaves. For type $j$, the expected error probability after querying model $m$ from belief $p$, and the reduction relative to immediate classification, are
\begin{align}
R_{jm}(p)&:=\sum_{x\in\{A,B\}}\min\{p\ell_{jm}^A(x),(1-p)\ell_{jm}^B(x)\} \mbox{ and } \Delta_{jm}^{\mathrm{1q}}(p):=\min\{p,1-p\}-R_{jm}(p) \label{eq:one-query-risk}
\end{align}
Each term in $R_{jm}(p)$ is the joint probability of output $x$ and a MAP error. The query has adjusted expected cost $\chi_m+b_jR_{jm}(p)$. Its net benefit relative to stopping is therefore $V_{jm}^{\mathrm{1q}}(p;\bm\chi):=b_j\Delta_{jm}^{\mathrm{1q}}(p)-\chi_m$. Proposition~\ref{prop:Bellman} then gives
\begin{equation*}
\Gamma_{j,1}(\bm\chi)=g_j(\rho_j)-\max\left\{0,\max_{m\in[M]}V_{jm}^{\mathrm{1q}}(\rho_j;\bm\chi)\right\}.
\end{equation*}
At any belief $p$, the job queries a model with the largest net benefit if that benefit is positive. If all net benefits are negative, it stops. At zero, stopping and the maximizing queries tie.

A simpler comparison emerges when accuracy is symmetric across labels: $q_{jm}^A=q_{jm}^B=q_{jm}>1/2$. Write $s(p):=\max\{p,1-p\}$ for the job's current MAP confidence. A model reduces expected classification error only if its accuracy exceeds this confidence. Accounting for its adjusted cost gives the index
\begin{equation*}
I_{jm}(\bm\chi):=q_{jm}-\frac{\chi_m}{b_j},
\qquad \kappa_j(\bm\chi):=\max_{m\in[M]}I_{jm}(\bm\chi).
\end{equation*}
\begin{proposition}[One-query index under symmetric accuracy]
\label{prop:L1-symmetric-index}
Suppose $L=1$ and $q_{jm}^A=q_{jm}^B=q_{jm}>1/2$. Then $\Delta_{jm}^{\mathrm{1q}}(p)=\bigl(q_{jm}-s(p)\bigr)^+$.
At fixed adjusted costs, stopping is optimal exactly when $s(p)\ge\kappa_j(\bm\chi)$. Querying model $m$ is optimal exactly when $I_{jm}(\bm\chi)=\kappa_j(\bm\chi)\ge s(p)$.
\end{proposition}

Under these assumptions, the posterior determines whether querying is worthwhile, but not which model is best when querying is optimal. That model maximizes accuracy minus the adjusted cost normalized by the error penalty. The ranking can still differ across job types and congestion prices. To obtain a system optimum, the resulting query decisions must generate the loads underlying those prices, as required by Corollary~\ref{cor:equilibrium-policy}.

When $b_j=b$ and $q_{jm}=q_m$ across types, the shared index makes it possible to solve directly for optimal loads. Types may still differ in their priors and arrival rates, and models may differ in their fees and capacities. Assigning more traffic to a model raises its congestion price and lowers its index, until the active models have a common adjusted index.

Let $s_j:=s(\rho_j)$ denote type $j$'s confidence at arrival and $\iota_m:=q_m-(c_m+h/\mu_m)/b$ denote model $m$'s index at zero load. For a proposed common index $\kappa$, an active model must satisfy $q_m-[c_m+\Theta_m(\Lambda_m)]/b=\kappa$. Solving this relation using \eqref{eq:Theta} gives
\begin{equation}
\Lambda_m(\kappa):=
\begin{cases}
\displaystyle\mu_m-\sqrt{\frac{h\mu_m}{b(q_m-\kappa)-c_m}},&\kappa<\iota_m,\\[2mm]
0,&\kappa\ge\iota_m.
\end{cases}
\label{eq:waterfill-load}
\end{equation}
A positive load can lower a model's index to $\kappa$ only if its idle index exceeds $\kappa$. Let $Q(\kappa):=\sum_m\Lambda_m(\kappa)$ be the aggregate query rate associated with these loads. The index rule also determines which types want to query. Their aggregate rates, excluding or including types indifferent between querying and stopping, are
\begin{equation*}
D_<(\kappa):=\sum_j\lambda_j\mathbf1\{s_j<\kappa\},
\qquad D_\le(\kappa):=\sum_j\lambda_j\mathbf1\{s_j\le\kappa\}.
\end{equation*}
Consistency requires $Q(\kappa)$ to lie between these two rates.

\begin{proposition}[Water-filling allocation]
\label{prop:L1-waterfill}
Suppose $L=1$, $b_j=b$, and $q_{jm}^A=q_{jm}^B=q_m>1/2$ for all $j,m$. The optimal loads are all zero if and only if $\max_m\iota_m\le\min_j s_j$. Otherwise, there is a unique scalar $\kappa^*\in[\min_j s_j,\max_m\iota_m)$ satisfying
\begin{equation}
D_<(\kappa^*)\le Q(\kappa^*)\le D_\le(\kappa^*),
\label{eq:waterfill-balance}
\end{equation}
and the unique optimal loads are $\Lambda_m^*=\Lambda_m(\kappa^*)$. Types with $s_j<\kappa^*$ query, types with $s_j>\kappa^*$ stop, and types with $s_j=\kappa^*$ may randomize between the two actions to satisfy \eqref{eq:waterfill-balance}. Every active model has adjusted index $\kappa^*$.
\end{proposition}

The scalar $\kappa^*$ determines both which types acquire information and how much traffic each model receives. Jobs with greater prior confidence stop, while those below the threshold query. On the model side, traffic lowers the indices of active models to the common level $\kappa^*$. Models whose idle indices do not exceed that level remain unused. This is the water-filling structure.

\begin{figure}[t]
\centering
\includegraphics[width=0.85\textwidth]{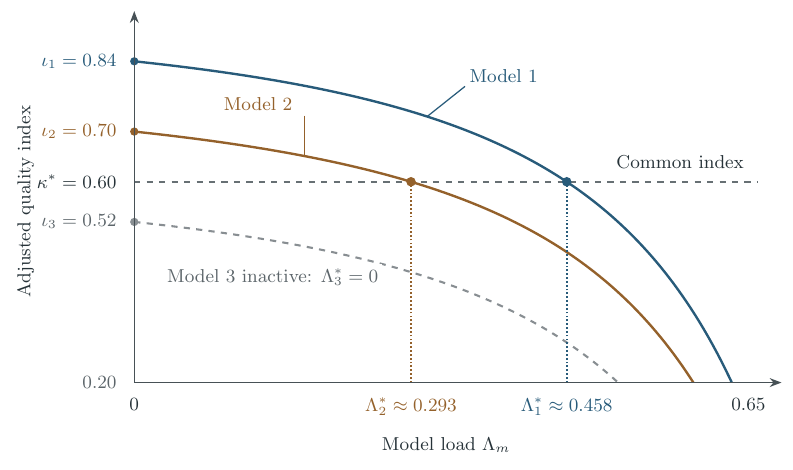}
\caption{One-query water-filling. Curves show $q_m-[c_m+\Theta_m(\Lambda_m)]/b$ for symmetric accuracies shared across types. Parameters: $b=10$, $h=1$, $\mu_m=1$, $(q_1,q_2,q_3)=(0.96,0.85,0.67)$, $(c_1,c_2,c_3)=(0.2,0.5,0.5)$, $(\lambda_1,\lambda_2)=(0.5,0.5)$, and $(s_1,s_2)=(0.5,0.6)$.}
\label{fig:water-filling}
\end{figure}

The inequalities in \eqref{eq:waterfill-balance} accommodate types whose confidence equals the common index. Only the fraction needed to match the model loads must query. In the positive-query case, one implementation sends each querying type to model $m$ with probability $\Lambda_m(\kappa^*)/Q(\kappa^*)$. Thus the scalar balance determines the optimal loads and also provides a routing policy that realizes them. Equation~\eqref{eq:waterfill-load} leaves a positive capacity buffer at every active model. Figure~\ref{fig:water-filling} illustrates the common index with two active models and a third whose idle index is too low to attract traffic.

The one-query comparison also determines the final query decision in a longer cascade, using the models still available. Appendix~\ref{subsec:M2-complete} applies it to the two-model, two-query case. 

\subsubsection*{Proofs for the one-query case}

\begin{proof}[Proof of Proposition~\ref{prop:L1-symmetric-index}]
Under symmetric accuracy, evaluating the two minima in
\eqref{eq:one-query-risk} gives
\[
R_{jm}(p)=\min\{p,1-p,1-q_{jm}\}.
\]
Subtracting from $\min\{p,1-p\}=1-s(p)$ yields
$\Delta_{jm}^{\mathrm{1q}}(p)=(q_{jm}-s(p))^+$.
If $q_{jm}\le s(p)$, querying $m$ has net benefit $-\chi_m<0$.
Otherwise, its net benefit is
$b_j[I_{jm}(\bm\chi)-s(p)]$. Thus stopping is optimal exactly when
every model's index is at most $s(p)$. A query is optimal exactly when
its index attains the maximum and is at least $s(p)$.
\end{proof}

\begin{proof}[Proof of Proposition~\ref{prop:L1-waterfill}]
At zero loads, model $m$ has adjusted cost $c_m+h/\mu_m$ and index $\iota_m$. By Proposition \ref{prop:L1-symmetric-index} and Corollary \ref{cor:equilibrium-policy}, immediate stopping for every type is system optimal exactly when $\max_m\iota_m\le\min_j s_j$. Uniqueness of optimal loads in Theorem \ref{thm:convexity} then gives the stated zero-query case.

Otherwise the attained optimum has positive aggregate query rate. Let $\kappa^*$ be the largest model index at its adjusted costs:
\[
\kappa^*:=\max_m
\left\{q_m-\frac{c_m+\Theta_m(\Lambda_m^*)}{b}\right\}.
\]
The index is common across types. Every queried model must attain this maximum by Proposition \ref{prop:L1-symmetric-index} and Corollary \ref{cor:equilibrium-policy}. For an active model, therefore,
\[
\Theta_m(\Lambda_m^*)=b(q_m-\kappa^*)-c_m.
\]
Since $\Theta_m$ is strictly increasing and $\Theta_m(0)=h/\mu_m$, the positive load exists exactly when $\kappa^*<\iota_m$, and inversion gives \eqref{eq:waterfill-load}. Every idle model has index
$\iota_m\le\kappa^*$. Thus $\Lambda_m^*=\Lambda_m(\kappa^*)$ for all $m$.

By the one-query index rule, types with $s_j<\kappa^*$ query, those
with $s_j>\kappa^*$ stop, and those at equality can mix. Summing query
rates proves \eqref{eq:waterfill-balance}. Positive flow implies
$\kappa^*\ge\min_j s_j$, and an active model implies
$\kappa^*<\max_m\iota_m$.

On this interval, $Q(\kappa)$ is positive, continuous, and strictly
decreasing: each active term in \eqref{eq:waterfill-load} decreases
strictly and vanishes continuously at its idle-index threshold.
If $\kappa_1<\kappa_2$ both satisfied \eqref{eq:waterfill-balance},
then
\[
Q(\kappa_1)\le D_\le(\kappa_1)
\le D_<(\kappa_2)\le Q(\kappa_2),
\]
contradicting strict decrease. This proves uniqueness.

For completeness, any scalar in the stated interval satisfying
\eqref{eq:waterfill-balance} implements a system optimum.
Choose type-wise query fractions equal to one when $s_j<\kappa$,
zero when $s_j>\kappa$, and fractions in $[0,1]$ at equality so that
their total query rate is $Q(\kappa)$. Route each querying type to
model $m$ with probability $\Lambda_m(\kappa)/Q(\kappa)$ and apply
MAP after the output. These choices induce exactly the specified
loads, which are strictly below capacity. All active models have
index $\kappa$, and idle models have index at most $\kappa$.
Every selected action is therefore optimal at the induced prices,
so Corollary \ref{cor:equilibrium-policy} establishes global optimality.
\end{proof}

\subsection{The two-model characterization}
\label{subsec:M2-complete}

For $M=L=2$, the organization chooses a first model and may conditionally
escalate to the other one. The one-query error probability $R_{jm}(p)$
and reduction $\Delta_{jm}^{\mathrm{1q}}(p)$ from Appendix~\ref{subsec:L1}
have the following explicit form when accuracies may differ across
the two labels.

\begin{lemma}[Exact one-query information value]
\label{lem:tent-shape}
Under Assumption \ref{ass:orientation}, let
\begin{equation*}
p_{jm}^{(1)}:=\frac{1-q_{jm}^B}{q_{jm}^A+1-q_{jm}^B},
\qquad
p_{jm}^{(2)}:=\frac{q_{jm}^B}{1-q_{jm}^A+q_{jm}^B}.
\end{equation*}
Then $p_{jm}^{(1)}<1/2<p_{jm}^{(2)}$ and
\begin{equation}
\Delta_{jm}^{\mathrm{1q}}(p)
=
\begin{cases}
0,
&0\le p\le p_{jm}^{(1)},
\\[1mm]
pq_{jm}^A-(1-p)(1-q_{jm}^B),
&p_{jm}^{(1)}<p\le1/2,
\\[1mm]
(1-p)q_{jm}^B-p(1-q_{jm}^A),
&1/2<p<p_{jm}^{(2)},
\\[1mm]
0,
&p_{jm}^{(2)}\le p\le1.
\end{cases}
\label{eq:tent-formula}
\end{equation}
Thus $\Delta_{jm}^{\mathrm{1q}}$ is a continuous tent-shaped function with
peak
\begin{equation}
\Delta_{jm}^{\mathrm{1q}}(1/2)
=\frac{q_{jm}^A+q_{jm}^B-1}{2}.
\label{eq:tent-peak}
\end{equation}
\end{lemma}

\begin{corollary}[Explicit one-query continuation interval]
\label{cor:one-query-cutoffs}
With only model $m$ available and one query remaining, querying is strictly
optimal exactly when
\begin{equation}
b_j\Delta_{jm}^{\mathrm{1q}}(p)>\chi_m.
\label{eq:one-query-rule}
\end{equation}
Let $\bar c_{jm}:=\chi_m/b_j$. If
$\bar c_{jm}\ge(q_{jm}^A+q_{jm}^B-1)/2$, querying is never strictly optimal. If the
inequality is reversed, querying is strictly optimal exactly for
\begin{equation}
\frac{\bar c_{jm}+1-q_{jm}^B}{q_{jm}^A+1-q_{jm}^B}
<p<
\frac{q_{jm}^B-\bar c_{jm}}{1-q_{jm}^A+q_{jm}^B}.
\label{eq:one-query-explicit-cutoffs}
\end{equation}
At equality in \eqref{eq:one-query-rule}, stopping and querying tie.
\end{corollary}

Thus a last query is valuable only over an intermediate posterior interval.
Increasing its adjusted cost weakly shrinks that interval.

For a candidate first model $i\in\{1,2\}$, write $k(i):=3-i$ and define
\begin{equation*}
p_{ji,x}:=\tau_{ji,x}(\rho_j),
\qquad
\bar s_{ji,x}:=s_{ji,x}(\rho_j),
\qquad x\in\{A,B\}.
\end{equation*}
The optimal cost conditional on querying model $i$ first is
\begin{equation}
F_{ji}^{(2)}(\bm\chi)
:=\chi_i+\sum_{x\in\{A,B\}}\bar s_{ji,x}
\min\left\{
g_j(p_{ji,x}),
\chi_{k(i)}+b_jR_{j,k(i)}(p_{ji,x})
\right\}.
\label{eq:M2-first-cost}
\end{equation}

\begin{theorem}[Complete two-model characterization]
\label{thm:M2-complete}
Fix type $j$, $M=L=2$, and adjusted costs $\bm\chi$. Then
\begin{equation}
\Gamma_{j,2}(\bm\chi)
=\min\left\{
g_j(\rho_j),F_{j1}^{(2)}(\bm\chi),F_{j2}^{(2)}(\bm\chi)
\right\}.
\label{eq:M2-root-value}
\end{equation}
If model $i$ is queried first, an optimal strict-continuation set is
\begin{equation}
C_{ji}^*(\bm\chi)
:=\left\{x\in\{A,B\}:
b_j\Delta_{j,k(i)}^{\mathrm{1q}}(p_{ji,x})>\chi_{k(i)}
\right\}.
\label{eq:M2-continuation-set}
\end{equation}
On a branch where equality holds, either stopping or querying the second model
is optimal. Consequently, up to irrelevant specifications below a stopping
node and a fixed tie convention for terminal labels, each type has nine routing
architectures: root stopping and, for each first model, escalation after
neither output, only after $A$, only after $B$, or after both. These
architectures can be evaluated in constant time for fixed $\bm\chi$.
\end{theorem}

The characterization separates the first-model comparison, which uses the
full expected continuation values in \eqref{eq:M2-first-cost}, from escalation,
which is governed by the posterior interval in
\eqref{eq:one-query-explicit-cutoffs}. Congestion affects both decisions through
the two adjusted costs.

\begin{corollary}[Equilibrium implementation with two models]
\label{cor:M2-equilibrium}
Suppose \(M=L=2\). At the unique optimal load vector, set
\(\chi_m^*=c_m+\Theta_m(\Lambda_m^*)\). Up to irrelevant actions below stopping nodes and terminal-label tie choices, every type-tree pair assigned positive mass has one of the nine routing architectures in Theorem \ref{thm:M2-complete}. Randomization, if needed to reproduce the optimal load vector, occurs only among architectures tied in adjusted cost.
\end{corollary}

For a two-model deployment, the remaining system decision is to choose traffic
shares among these explicit architectures while matching the adjusted costs
to the loads they induce.

\begin{proof}[Proof of Lemma~\ref{lem:tent-shape}]
The two minima in \eqref{eq:one-query-risk} change branches at the solutions of
\[
pq_{jm}^A=(1-p)(1-q_{jm}^B),
\qquad
p(1-q_{jm}^A)=(1-p)q_{jm}^B,
\]
which are $p_{jm}^{(1)}$ and $p_{jm}^{(2)}$, respectively. Positive
orientation is equivalent to
$p_{jm}^{(1)}<1/2<p_{jm}^{(2)}$. Evaluating the minima on the four resulting
intervals gives
\begin{align*}
R_{jm}(p)
&=p,
&&0\le p\le p_{jm}^{(1)},\\
&=(1-p)(1-q_{jm}^B)+p(1-q_{jm}^A),
&&p_{jm}^{(1)}<p<p_{jm}^{(2)},\\
&=1-p,
&&p_{jm}^{(2)}\le p\le1.
\end{align*}
Subtracting these expressions from $\min\{p,1-p\}$ proves
\eqref{eq:tent-formula}. The middle-left expression has slope
$q_{jm}^A+1-q_{jm}^B>0$, and the middle-right expression has slope
$-(1-q_{jm}^A+q_{jm}^B)<0$. The pieces agree at the breakpoints, and
evaluation at $p=1/2$ gives \eqref{eq:tent-peak}.
\end{proof}

\begin{proof}[Proof of Corollary~\ref{cor:one-query-cutoffs}]
The cost difference between stopping and querying is
\[
b_j\min\{p,1-p\}-\{\chi_m+b_jR_{jm}(p)\}
=b_j\Delta_{jm}^{\mathrm{1q}}(p)-\chi_m.
\]
This proves \eqref{eq:one-query-rule}. Comparing
$\bar c_{jm}=\chi_m/b_j$ with the
peak in \eqref{eq:tent-peak}, and solving the two affine equalities in
\eqref{eq:tent-formula}, gives \eqref{eq:one-query-explicit-cutoffs} and the
stated tie convention.
\end{proof}

\begin{proof}[Proof of Theorem~\ref{thm:M2-complete}]
Suppose model $i$ is queried first. Conditional on output $x$, the posterior is
$p_{ji,x}$, the only remaining model is $k(i)$, and one query remains. By
Corollary \ref{cor:one-query-cutoffs}, the branch value is
\[
\min\left\{g_j(p_{ji,x}),
\chi_{k(i)}+b_jR_{j,k(i)}(p_{ji,x})\right\},
\]
and strict continuation occurs exactly on the set in
\eqref{eq:M2-continuation-set}. Taking expectation with weights
$\bar s_{ji,x}$ and adding $\chi_i$ gives $F_{ji}^{(2)}$. Minimizing over root
stopping and $i\in\{1,2\}$ proves \eqref{eq:M2-root-value}. For each first
model, the two output branches can independently stop or continue, yielding
four continuation sets and hence $1+2\cdot4=9$ routing architectures. Actions
below stopping nodes and terminal-label choices at MAP ties do not affect this
architectural count.
\end{proof}

\begin{proof}[Proof of Corollary~\ref{cor:M2-equilibrium}]
At the unique optimal load vector, Corollary \ref{cor:equilibrium-policy} evaluates the
type-wise problem at adjusted costs
\(\chi_m^*=c_m+\Theta_m(\Lambda_m^*)\), \(m\in[M]\).
Theorem \ref{thm:M2-complete} exhausts their routing architectures. The congestion-price conditions in Theorem \ref{thm:tree-KKT} therefore assign positive
mass only to cost-minimizing architectures whose mixture induces
$\bm\Lambda^*$.
\end{proof}

\FloatBarrier
\section{Extension: Repeated Queries with Queue Reentry}
\label{sec:repeated-queries}
\label{apdx:repeated-queries}

So far our main problem lets a job consult each model at most once. We now extend it to allow returning to a model for another output. This expands the routing choices, but it also requires a description of how repeated outputs are related. The queueing and convex optimization results need only a specified joint output law. The simpler Bayesian state requires the stronger condition that each call provides a fresh signal conditional on the true label. We separate these two requirements below.

\subsection{Repeated calls and system costs}

Fix a finite total query cap $L\ge1$, which may exceed $M$. A job may stop without making a query. After each completed call, it releases the server. If it queries the same model again, the new call joins the tail of that model's FCFS queue and receives a new independent exponential service requirement with rate $\mu_m$. Thus the job retains neither the server nor its place in the queue between calls. For example, a realized route $(3,7,3,3,5)$ makes five visits, three of them to model $3$. The physical model is revisited, but each completion advances the job to a new stage of its finite route. Note that every new call to the same model must rejoin the queue. Consecutive repeats may receive immediate service if the queue is empty. The results we discuss below do not cover a job that keeps the server to obtain further outputs while other jobs wait.

Let $X_{jm}^{(k)}\in\{A,B\}$ denote the potential output of the $k$th call to model $m$ by a type-$j$ job. Conditional on $Y_j$, the finite array $(X_{jm}^{(k)}:m\in[M],k\in[L])$ has a specified joint law. Arrays are independent across jobs and independent of arrivals, service requirements, queue states, and policy randomization. At this stage, repeated outputs may be correlated.

Let $\cT_L^{\mathrm{rep}}$ be the finite set of deterministic CAAR trees with at most $L$ queries on any path. These trees may stop immediately, and every model remains available after every query. For a tree $T$, let $N_{jTm}$ count its visits to model $m$. Its error probability $e_{jT}$ is defined as before, while its visit coefficient becomes $v_{jTm}:=\E[N_{jTm}]$. This expected count can exceed one. A tree mixture $\bm z$ induces the load
\[
\Lambda_m^{\mathrm{rep}}(\bm z)
=\sum_j\lambda_j\sum_{T\in\cT_L^{\mathrm{rep}}}z_{jT}v_{jTm}.
\]
Recording each job's route and current stage again gives an open BCMP network, as in Section~\ref{sec:system-foundation}. The exact holding-cost formula therefore continues to apply. With $\omega_{jT}:=b_je_{jT}+\sum_m c_mv_{jTm}$, the objective retains the form
\[
C_L^{\mathrm{rep}}(\bm z)
=\sum_j\lambda_j\sum_{T\in\cT_L^{\mathrm{rep}}}z_{jT}\omega_{jT}
+\sum_m\Psi_m\bigl(\Lambda_m^{\mathrm{rep}}(\bm z)\bigr).
\]
As in Section~\ref{sec:system-foundation}, mixtures satisfy the type-wise simplex constraints, and $\Psi_m$ is infinite at unstable loads.

\begin{proposition}[System structure under repeated queries]
\label{prop:repeated-query-extension}
Under the queue-reentry convention and the specified output law, the following statements hold.
\begin{enumerate}[label=(\roman*),leftmargin=*]
\item Every randomized repeated-query CAAR policy is equivalent to a mixture over $\cT_L^{\mathrm{rep}}$. If $\Lambda_m^{\mathrm{rep}}(\bm z)<\mu_m$ for every $m$, its exact holding-cost rate is $\sum_m\Psi_m(\Lambda_m^{\mathrm{rep}}(\bm z))$.
\item The objective $C_L^{\mathrm{rep}}$ is convex and attains its minimum. All optimal mixtures induce the same physical load vector. A stable mixture $\bm z^*$ is optimal if and only if every positively used type-$j$ tree minimizes over $T\in\cT_L^{\mathrm{rep}}$:
\[
A_{jT}^{\mathrm{rep}}(\bm\theta^*)
:=b_je_{jT}+\sum_m(c_m+\theta_m^*)v_{jTm},
\qquad
\theta_m^*=\Theta_m\bigl(\Lambda_m^{\mathrm{rep}}(\bm z^*)\bigr).
\]
\item There exists an optimal mixture with at most $M+J$ positive type--tree weights.
\end{enumerate}
\end{proposition}

The change from visit probabilities to visit counts is what incorporates repeated calls into both fees and congestion. Each repeat pays the same marginal congestion price as another visit to that physical model. The queueing argument requires neither independent outputs within a job nor Poisson internal arrival streams. Immediate stopping also preserves a finite-cost feasible policy at any external arrival rates. The support bound from Proposition~\ref{thm:sparsity} continues to hold because loads remain linear in the tree probabilities.

\begin{proof}[Proof of Proposition~\ref{prop:repeated-query-extension}]
Fix a type. Pre-sample the policy's entire private random tape at arrival. Conditional on that tape, the action at every finite output history is fixed, so the policy determines one complete tree in $\cT_L^{\mathrm{rep}}$. Mixing over the tape gives the tree-mixture representation. Conversely, sampling a tree at arrival implements any prescribed mixture.

For the queueing representation, also pre-sample the label and the finite potential-output array. This arrival-time mark determines a realized route $\mathfrak r=(m_1,\ldots,m_d)$ with $0\le d\le L$, where model indices may repeat. Retain the positive-probability nonempty routes. Independent Poisson marking splits the external arrivals into route streams. A customer class records its type, route, and current route stage. Service at stage $t$ is exponential with rate $\mu_{m_t}$, and completion moves the customer to stage $t+1$ or out of the network. If $m_{t+1}=m_t$, the class changes and the job returns to the tail of the same physical queue. The route-stage index strictly increases, so the reachable class graph is finite and open. Every physical node is a type-1 FCFS BCMP station because all classes at model $m$ have the same exponential rate $\mu_m$.

The sum of the route-stage throughputs at model $m$ is exactly $\Lambda_m^{\mathrm{rep}}(\bm z)$ because every appearance of $m$ on a route contributes one visit. Summing the local BCMP factors over ordered class sequences of length $n$ gives $(\Lambda_m^{\mathrm{rep}}(\bm z)/\mu_m)^n$. Hence the local normalizing series is finite under the stated stability inequalities. If $\mathcal N_m$ is the stationary number of jobs at model $m$, then
\[
\E[\mathcal N_m]
=\frac{\Lambda_m^{\mathrm{rep}}(\bm z)}
{\mu_m-\Lambda_m^{\mathrm{rep}}(\bm z)}.
\]
Conversely, stationarity would imply $\Lambda_m^{\mathrm{rep}}=\mu_m\Pp(\mathcal N_m>0)<\mu_m$ at every active model. Strictness follows because the empty state belongs to the irreducible reachable class and has positive stationary probability. Multiplication by $h$ and summation over models proves part (i), including the instability convention in the extended-valued objective.

For part (ii), the objective $C_L^{\mathrm{rep}}$ is defined on a finite product of simplices. The immediate-stopping trees give a finite-cost feasible point. Each load is affine in $\bm z$, and the noncongestion cost is linear. Lower semicontinuity and convexity of the functions $\Psi_m$ therefore give an attained minimum. If two optima had different load vectors, their midpoint would have the average linear cost and a strictly smaller holding cost, because each $\Psi_m$ is strictly convex on its stable domain. Thus optimal loads are unique.

Write $\Lambda_m^*:=\Lambda_m^{\mathrm{rep}}(\bm z^*)$. At an optimal stable mixture, shifting a small amount of probability from an active tree $T$ to another tree $T'$ of the same type gives
\[
\lambda_j\bigl[A_{jT'}^{\mathrm{rep}}(\bm\theta^*)
-A_{jT}^{\mathrm{rep}}(\bm\theta^*)\bigr]\ge0.
\]
This is the right directional derivative of the objective. If an optimal load is zero, every active tree has zero visits to that model, so the required derivative is the right derivative $\Theta_m(0)$. Necessity of the minimizing condition follows. Conversely, the supporting inequality
\[
\Psi_m(\Lambda)-\Psi_m(\Lambda_m^*)
\ge\Theta_m(\Lambda_m^*)(\Lambda-\Lambda_m^*)
\]
and the type-wise minimizing condition show that no feasible mixture has smaller cost. This proves part (ii).

For part (iii), apply the fixed-load argument in the proof of Proposition~\ref{thm:sparsity} to an optimal repeated-query mixture. Its tree coefficients are expected visit counts under the specified joint output law. They remain fixed when the mixture weights change, so preserving the $M$ model loads and $J$ probability sums again gives an extreme point with at most $M+J$ positive weights. This argument does not require independent repeated outputs.
\end{proof}

\subsection{Routing and stopping with repeated observations}
\label{sec:repeated-routing}

To recover a compact Bayesian recursion, suppose now that, conditional on $Y_j$, all coordinates of the potential-output array are mutually independent and
\[
\Pp(X_{jm}^{(k)}=x\mid Y_j=y)=\ell_{jm}^y(x)
\quad\text{for every }k.
\]
The likelihood does not change with the replication count. This is an additional assumption about repeated calls. It does not follow from the baseline model, which specifies only one output per model. Proposition~\ref{prop:repeated-Bayesian} below formalizes the resulting Bayesian and stopping properties.

Under this condition, the posterior $p$ determines the predictive law of every next call through $s_{jm,x}(p)$ and $\tau_{jm,x}(p)$. Every model remains available. There is therefore no unused-model set $S$ to carry forward, and neither the most recently queried model nor the replication counts are needed. The continuation state is $(p,r)$, where $r$ is the remaining query budget.

For fixed adjusted costs $\bm\chi\in(0,\infty)^M$, let $W_{j,r}^{\mathrm{rep}}(p;\bm\chi)$ be the minimum remaining adjusted cost. The terminal value is $W_{j,0}^{\mathrm{rep}}(p;\bm\chi)=g_j(p)$, and for $r\ge1$ the Bellman recursion is
\begin{equation}
W_{j,r}^{\mathrm{rep}}(p;\bm\chi)
=\min\left\{
 g_j(p),\ 
 \min_{m\in[M]}\left[
 \chi_m+\sum_{x\in\{A,B\}}s_{jm,x}(p)
 W_{j,r-1}^{\mathrm{rep}}(\tau_{jm,x}(p);\bm\chi)
 \right]
\right\}.
\label{eq:repeated-query-Bellman}
\end{equation}
The root value satisfies
\[
W_{j,L}^{\mathrm{rep}}(\rho_j;\bm\chi)
=\min_{T\in\cT_L^{\mathrm{rep}}}
\left\{b_je_{jT}+\sum_m\chi_mv_{jTm}\right\}.
\]
Stopping remains an option at the root. After a query, the policy compares stopping with all models again, including the model just queried. At the optimal system prices, these decisions characterize the trees used in Proposition~\ref{prop:repeated-query-extension}.

The stopping structure from Section~\ref{sec:policy-structure} also carries over. For each fixed $j,r,\bm\chi$, there are thresholds $0\le p^-\le1/2\le p^+\le1$ such that stopping is optimal exactly on $[0,p^-]\cup[p^+,1]$. If $\bm\chi'\ge\bm\chi$ componentwise, the region where continuation is strictly better than stopping under $\bm\chi'$ is contained in the corresponding region under $\bm\chi$. Higher adjusted query costs therefore weakly reduce the range of posteriors at which another call, including a repeat, is worthwhile. This comparison holds for the fixed-price continuation problem.

\begin{proposition}[Bayesian decisions with repeated queries]
\label{prop:repeated-Bayesian}
Under the stationary conditional independence assumption in Appendix~\ref{sec:repeated-routing}, the following statements hold for fixed type $j$ and adjusted costs $\bm\chi\in(0,\infty)^M$.
\begin{enumerate}[label=(\roman*),leftmargin=*]
\item The state $(p,r)$ is sufficient for optimal continuation. The value satisfies \eqref{eq:repeated-query-Bellman}, with $W_{j,0}^{\mathrm{rep}}=g_j$. Its root value $W_{j,L}^{\mathrm{rep}}(\rho_j;\bm\chi)$ equals the minimum adjusted cost over $\cT_L^{\mathrm{rep}}$.
\item For each $r$, there are thresholds $0\le p^-\le1/2\le p^+\le1$ such that stopping is optimal exactly on $[0,p^-]\cup[p^+,1]$. If $\bm\chi'\ge\bm\chi$ componentwise, then
\[
\{p:W_{j,r}^{\mathrm{rep}}(p;\bm\chi')<g_j(p)\}
\subseteq
\{p:W_{j,r}^{\mathrm{rep}}(p;\bm\chi)<g_j(p)\}.
\]
\end{enumerate}
\end{proposition}

\begin{proof}[Proof of Proposition~\ref{prop:repeated-Bayesian}]
Consider a reachable history $\mathfrak h_t=((m_1,x_1),\ldots,(m_t,x_t))$ and let $\kappa_s:=\sum_{u=1}^s\mathbf1\{m_u=m_s\}$. Its stage-$s$ output is the coordinate $X_{jm_s}^{(\kappa_s)}$. The pairs $(m_s,\kappa_s)$ are distinct. Write the stage-$s$ action as $\varphi_s(u,\mathfrak h_{s-1})$. For any Borel set $B$ of private seeds, define
\[
a_B(\mathfrak h_t)
:=\int_B\prod_{s=1}^t
\mathbf1\{\varphi_s(u,\mathfrak h_{s-1})=m_s\}\,du.
\]
The joint probability of the label, seed event, and history is
\[
\Pp(Y_j=y)a_B(\mathfrak h_t)
\prod_{s=1}^t\ell_{jm_s}^y(x_s).
\]
Because this identity holds for every $B$, the label-independent selection term cancels in the posterior conditional on the full private seed and history. Conditional independence also leaves every unrevealed coordinate with its original likelihood. Therefore the scalar posterior determines the predictive law of a fresh call to every model, regardless of how often that model was previously queried.

At any stopping history, replacing the terminal decision by MAP weakly lowers error and does not alter visits, so the stopping loss is $g_j(p)$. Backward induction on $r$ proves \eqref{eq:repeated-query-Bellman}. The first action either stops at cost $g_j(p)$ or queries model $m$, pays $\chi_m$, observes $x$ with probability $s_{jm,x}(p)$, and continues with $r-1$ calls. The same induction identifies the root value with the minimum adjusted cost over deterministic repeated-query trees. Immediate stopping is permitted in both the recursion and the tree set.

For the stopping thresholds, fix $j,r,\bm\chi$. There are finitely many deterministic continuation trees. Conditional on $Y_j=A$ and on $Y_j=B$, a fixed tree has two fixed expected costs. Its Bayes cost is therefore affine in $p$. Consequently, $W_{j,r}^{\mathrm{rep}}(p;\bm\chi)$ is the pointwise minimum of finitely many affine functions, and is continuous and concave. Immediate MAP stopping gives
\[
0\le W_{j,r}^{\mathrm{rep}}(p;\bm\chi)\le g_j(p),
\qquad
W_{j,r}^{\mathrm{rep}}(0;\bm\chi)
=W_{j,r}^{\mathrm{rep}}(1;\bm\chi)=0.
\]
On $[0,1/2]$, the function $b_jp-W_{j,r}^{\mathrm{rep}}(p;\bm\chi)$ is nonnegative and convex, and vanishes at zero. Its zero set is therefore a closed interval $[0,p^-]$. On $[1/2,1]$, the same argument applied to $b_j(1-p)-W_{j,r}^{\mathrm{rep}}(p;\bm\chi)$ gives a closed interval $[p^+,1]$. Stopping is optimal exactly on the union of these intervals.

Finally, increasing $\bm\chi$ componentwise weakly increases the expected cost of every fixed continuation policy and leaves the stopping loss unchanged. Taking the minimum over policies gives
\[
W_{j,r}^{\mathrm{rep}}(p;\bm\chi')
\ge W_{j,r}^{\mathrm{rep}}(p;\bm\chi)
\quad\text{when }\bm\chi'\ge\bm\chi.
\]
If continuation is strictly better than stopping under $\bm\chi'$, it is also strictly better under $\bm\chi$. This proves the stated inclusion of continuation regions.
\end{proof}

\subsection{Information screening and load bounds}

Under Assumption~\ref{ass:orientation}, the information screen in \eqref{eq:statewise-frontier} and the model exclusions and load bound in Corollary~\ref{cor:exogenous-frontier} also extend to repeated calls. Here the model load counts every visit, including repeats. Repeated query therefore preserves the capacity protection created by its limited information value. Proposition~\ref{prop:repeated-screening} below states this extension.

The information assumption is essential for the reduced state and these screening conclusions. If the law of the next output depends on previous calls beyond their information about the true label, the policy must retain a sufficient statistic for that dependence. The queueing and convexity results still apply under the specified joint output law above.

\begin{proposition}[Information limits on repeated queries]
\label{prop:repeated-screening}
Assume the stationary conditional independence condition in Section~\ref{sec:repeated-routing} and Assumption~\ref{ass:orientation}.
\begin{enumerate}[label=(\roman*),leftmargin=*]
\item If querying model $m$ is Bellman optimal at $(p,r)$ with $r\ge1$ and adjusted costs $\bm\chi\in(0,\infty)^M$, then
\[
\chi_m\le b_j\min\{p,1-p,2p(1-p)d_{jm}\}.
\]
\item For any optimal repeated-query mixture $\bm z^*$, write $\Lambda_m^*:=\Lambda_m^{\mathrm{rep}}(\bm z^*)$. With $\cF_j^0$, $\cF^0$, and $\overline v_m$ defined in Section~\ref{subsec:information-cost-frontier},
\[
m\notin\cF_j^0
\quad\Longrightarrow\quad
\sum_{T\in\cT_L^{\mathrm{rep}}}z_{jT}^*v_{jTm}=0.
\]
Restricting each type's trees to models in $\cF_j^0$ preserves the optimal value. If $m\notin\cF^0$, then $\Lambda_m^*=0$. Every active model satisfies $\overline v_m>\chi_m^{\mathrm{idle}}$ and
\[
\Lambda_m^*\le\mu_m-\sqrt{\frac{h\mu_m}{\overline v_m-c_m}}.
\]
\end{enumerate}
\end{proposition}

\begin{proof}
Fix $j,r,\bm\chi$ with $r\ge1$, and write $V(p):=W_{j,r-1}^{\mathrm{rep}}(p;\bm\chi)$. The proof of Proposition~\ref{prop:repeated-Bayesian} shows that $V$ is continuous and concave, with $V(0)=V(1)=0$ and $0\le V(p)\le b_j\min\{p,1-p\}$. For $0<u<v<1$, concavity and these bounds give
\[
-b_j\le-\frac{V(v)}{1-v}
\le\frac{V(v)-V(u)}{v-u}
\le\frac{V(u)}{u}\le b_j.
\]
Continuity covers the endpoints. Thus $|V(p)-V(p')|\le b_j|p-p'|$ for all $p,p'\in[0,1]$.

Let $P'$ be the posterior after a query to $m$ from belief $p$. The predictive law and posterior update are the same as in Section~\ref{subsec:information-cost-frontier}, so
\[
\E[|P'-p|]
=\sum_{x\in\{A,B\}}s_{jm,x}(p)|\tau_{jm,x}(p)-p|
=2p(1-p)d_{jm}.
\]
A policy allowed at most $r$ calls can instead use an optimal policy with at most $r-1$ calls. Every model remains available, and fresh future calls have the same likelihoods under either choice. Hence $W_{j,r}^{\mathrm{rep}}(p;\bm\chi)\le V(p)$. If querying $m$ is optimal, including a tie, the Bellman recursion gives
\[
\chi_m+\E[V(P')]
=W_{j,r}^{\mathrm{rep}}(p;\bm\chi)
\le V(p).
\]
It follows that
\[
\chi_m\le V(p)-\E[V(P')]
\le\E[|V(p)-V(P')|]
\le2b_jp(1-p)d_{jm}.
\]
Immediate stopping also gives $\chi_m\le W_{j,r}^{\mathrm{rep}}(p;\bm\chi)\le g_j(p)$ because future costs are nonnegative. Combining these inequalities proves part (i).

For part (ii), let $\chi_m^*:=c_m+\Theta_m(\Lambda_m^*)$. Proposition~\ref{prop:repeated-query-extension} implies that each tree with positive mixture weight minimizes adjusted cost at $\bm\chi^*$. Every continuation at a history reached with positive probability by that tree must therefore be Bellman optimal. Otherwise, replacing the continuation by a strictly cheaper one would strictly reduce the tree's expected adjusted cost. This replacement is feasible because repeated queries leave every model available and the remaining budget is unchanged.

If type $j$ makes a positive expected number of visits to $m$, some positively weighted tree queries $m$ at a history it reaches with positive probability. Part (i) and $2p(1-p)\le1/2$ imply $\chi_m^*\le b_jd_{jm}/2$. Since $\lambda_j>0$, those visits imply $\Lambda_m^*>0$ and therefore
\[
\chi_m^{\mathrm{idle}}
=c_m+\Theta_m(0)
<\chi_m^*\le\frac{b_jd_{jm}}2.
\]
Thus positive use requires $m\in\cF_j^0$. In particular, equality between the idle cost and the information bound is insufficient for positive traffic, even if some fixed-price decisions admit ties. Models outside $\cF_j^0$ have zero expected visits for type $j$. Any occurrences of them at histories of probability zero can be replaced by stopping without changing costs or loads. The resulting trees use only models in $\cF_j^0$, which proves preservation of the optimal value under this restriction. Summing over types proves that $m\notin\cF^0$ implies $\Lambda_m^*=0$.

Finally, if $\Lambda_m^*>0$, at least one type uses $m$ with positive expected visit count. The preceding argument gives
\[
\chi_m^{\mathrm{idle}}<\chi_m^*
=c_m+\frac{h\mu_m}{(\mu_m-\Lambda_m^*)^2}
\le\overline v_m.
\]
The denominator $\overline v_m-c_m$ is positive. Rearranging gives the claimed load bound. This argument counts every repeated visit and does not require $v_{jTm}\le1$.
\end{proof}

\end{document}